\documentclass{agtart_a}
\pdfoutput=1

\usepackage[percent]{overpic}

\title[Regular homotopy and total curvature I]
{Regular homotopy and total curvature I:\\ circle immersions into surfaces}

\author{Tobias Ekholm}
\givenname{Tobias}
\surname{Ekholm}
\address{Department of mathematics\\USC\\\newline
Los Angeles CA 90803\\
USA}
\email{tekholm@usc.edu}

\volumenumber{6}
\issuenumber{}
\publicationyear{2006}
\papernumber{16}
\startpage{459}
\endpage{492}

\doi{}
\MR{}
\Zbl{}

\keyword{circle immersion}
\keyword{geodesic curvature}
\keyword{regular curve}
\keyword{regular homotopy}
\keyword{Riemann surface}
\keyword{total curvature}
\subject{primary}{msc2000}{53C42}
\subject{secondary}{msc2000}{53A04}
\subject{secondary}{msc2000}{57R42}

\received{8 February 2005}
\revised{22 February 2006}
\accepted{12 March 2006}
\published{23 March 2006}
\publishedonline{23 March 2006}
\proposed{}
\seconded{}
\corresponding{}
\editor{}
\version{}

\arxivreference{math.GT/0310266}

\AtBeginDocument{\let\bar\wbar\let\tilde\wtilde\let\hat\what}

\makeatletter
\def\cnewtheorem#1[#2]#3{\newtheorem{#1}{#3}[section]
\expandafter\let\csname c@#1\endcsname\c@thm}

\renewcommand{\baselinestretch}{1.04}  

\newtheorem{thm}{Theorem}[section]
\cnewtheorem{prp}[thm]{Proposition}
\cnewtheorem{cor}[thm]{Corollary}
\cnewtheorem{lma}[thm]{Lemma}

\theoremstyle{definition}
\cnewtheorem{dfn}[thm]{Definition}
\cnewtheorem{rmk}[thm]{Remark}

\makeatother  

\makeautorefname{prp}{Proposition}
\makeautorefname{lma}{Lemma}

\numberwithin{equation}{section}
\numberwithin{section}{part}
\renewcommand{\thesection}{\arabic{section}}

\renewcommand{\S}{{\mathbb{S}}}
\newcommand{\la}{\langle}
\newcommand{\ra}{\rangle}
\newcommand{\pa}{\partial}
\newcommand{\id}{\operatorname{id}}
\newcommand{\krn}{\operatorname{ker}}

\newcommand{\inr}{\operatorname{int}}

\newcommand{\rank}{\operatorname{rank}}
\newcommand{\vol}{\operatorname{vol}}
\newcommand{\pgc}{{\textsc{pgc}}}
\newcommand{\area}{\operatorname{area}}
\newcommand{\Ordo}{\mathcal{O}}

\newcommand{\Imm}{\operatorname{Imm}}
\newcommand{\pImm}{\operatorname{Pgc}}

\begin{document}

\part{Circle immersions into surfaces}

\begin{abstract}
We consider properties of the total absolute geodesic curvature
functional on circle immersions into a Riemann surface. In particular,
we study its behavior under regular homotopies, its infima in regular
homotopy classes, and the homotopy types of spaces of its local
minima.
\end{abstract}

\maketitle

\section{Introduction}\label{sec:1}
An {\em immersion} of manifolds is a map with everywhere injective
differential. Two immersions are {\em regularly homotopic} if there
exists a continuous $1$--parameter family of immersions connecting one
to the other. The Smale--Hirsch $h$--principle \cite{Sm2,Hi} says
that the space of immersions $M\to N$, $\dim(M)<\dim(N)$ is homotopy
equivalent to the space of injective bundle maps $TM\to TN$. In
contrast to differential topological properties, differential
geometric properties of immersions do not in general satisfy
$h$--principles, see \cite[(A) on page 62]{Gr}. In this paper and the
sequel \cite{EII}, we study some aspects of the differential geometry
of immersions and regular homotopies in the most basic cases of
codimension one immersions. We investigate whether or not it is
possible to perform topological constructions while keeping control of
certain geometric quantities.

Let $\Sigma$ be a Riemann surface, ie, an orientable $2$--manifold with a Riemannian metric, and let  $c\co S^1\to\Sigma$ be an immersion of the circle parameterized by arc length. If $\vec c\co S^1\to U\Sigma$, where $U\Sigma$ is the unit tangent bundle of $\Sigma$, denotes the natural lift of $c$, then the $h$--principle mentioned above implies that the map $c\mapsto \vec c$ induces a weak homotopy equivalence between the space of circle immersions into $\Sigma$ and the space of continuous circle maps into $U\Sigma$. In particular, regular homotopy classes of circle immersions into $\Sigma$ are in one to one correspondence with the homotopy classes of (free) loops in $U\Sigma$.

The {\em total absolute geodesic curvature} $\kappa$ of a circle immersion $c$ into a Riemann surface is given by the integral
\begin{equation*}
\kappa(c)=\int_c |k_g|\,ds,
\end{equation*}
where $k_g$ is the geodesic curvature of $c$, and where $ds$ denotes the arc length element along $c$. We study properties of the functional $\kappa$ on the space of circle immersions, starting with the following question. {\em If $c_0$ and $c_1$ are regularly homotopic circle immersions into a Riemann surface, what is the infimum, over all regular homotopies $c_t$, $0\le t\le 1$, connecting $c_0$ to $c_1$, of $\max_{0\le t\le 1}\kappa(c_t)$?}

\fullref{thmK=const} answers this question for the simplest Riemann surfaces of constant curvature. We use the following notational conventions: All Riemann surfaces are assumed to be complete unless otherwise explicitly stated. If $\Sigma$ is a Riemann surface then $K\co\Sigma\to\R$ denotes its Gaussian curvature function. For topological spaces $X$ and $Y$, we write $X\approx Y$ to indicate that $X$ is homeomorphic to $Y$.
\begin{thm}\label{thmK=const}
Let $\Sigma$ be a Riemann surface and let $c_0,c_1\co S^1\to\Sigma$
be regularly homotopic.
\begin{itemize}
\item[{\rm (a)}]
If $\Sigma$ has constant curvature $K=0$ (the completeness assumption then implies $\Sigma\approx \R^2$, $\Sigma\approx S^1\times\R$, or $\Sigma\approx T^2$), or if $\Sigma\approx \R^2$ and has constant curvature $K<0$, then there exists a regular homotopy $c_t$, $0\le t\le 1$,
connecting $c_0$ to $c_1$ with
\begin{equation*}
\kappa(c_t)\le \max\{\kappa(c_0),\kappa(c_1)\},\quad 0\le t\le 1.
\end{equation*}
\item[{\rm (b)}]
If $\Sigma$ is the $2$--sphere with a constant curvature metric (with $K>0$) then, for any $\epsilon>0$, there exists a regular homotopy $c_t$, $0\le t\le 1$, connecting $c_0$ to $c_1$ such that
\begin{equation*}\label{eq1S^2}
\kappa(c_t)\le\max\{\kappa(c_0),\kappa(c_1),2\pi+\epsilon\},\quad
0\le t\le 1.
\end{equation*}
Moreover, if $c_0\co S^1\to\Sigma$ runs $m$ times around a geodesic and $c_1\co S^1\to\Sigma$ runs $m+2$ times around a geodesic then any regular homotopy $c_t$, $0\le t\le 1$, connecting $c_0$ to $c_1$ has an instant $c_{\tau}$, $0<\tau<1$, with
\begin{equation*}\label{eq2S^2}
\kappa(c_\tau)>2\pi.
\end{equation*}
\end{itemize}
\end{thm}
\fullref{thmK=const} is proved in \fullref{ssec:5C}. The
proof of (b) uses Arnol'd's $J^{-}$--invariant for immersed curves on
the sphere, see Arnol'd \cite{Ar}, Inshakov \cite{In} and Tchernov
\cite{Tch}. In \fullref{rmkKle0} we present a metric on $\R^2$ with
$K\le 0$ for which the conclusion in (a) does not hold.

The proof of \fullref{thmK=const} also gives information about infima of $\kappa$. To state these results we first introduce some notation. If $\Sigma$ is a flat Riemann surface then parallel translation gives a trivialization of $U\Sigma$ and the free homotopy classes of curves in $U\Sigma$ are in natural one to one correspondence with $\pi_1(\Sigma)\times\Z$, where $\pi_1(\Sigma)$ encodes the homotopy class of a circle immersion and $\Z$ its tangential degree. Thus, if $\Sigma\approx\R^2$, then we denote a regular homotopy class of circle immersions by the integer $m$ which equals the tangential degree of any of its representatives, and, similarly, if $\Sigma\approx S^1\times\R$ or $\Sigma\approx T^2$, then we denote a regular homotopy class by $(\xi,m)\in\pi_1(\Sigma)\times\Z$, where $\xi$ and $m$ is the homotopy class in $\Sigma$ and the tangential degree, respectively, of any of its representatives. If $\Sigma$ is the $2$--sphere then there are exactly two regular homotopy classes: one represented by a simple closed curve, the other by such a curve traversed twice. Finally, if $\alpha$ is a regular homotopy class of curves in a Riemann surface then let $\hat\kappa(\alpha)=\inf_{c\in\alpha}\kappa(c)$.

\begin{thm}\label{thminf}

\noindent
\begin{itemize}
\item[{\rm (a)}]
Let $\Sigma$ be a Riemann surface with $K(p)<0$ for all $p\in\Sigma$ and assume that either $\Sigma$ is closed or $\Sigma\approx\R^2$. Then the infimum $\hat\kappa(\alpha)$ is attained at some curve in the regular homotopy class $\alpha$ if and only if $\alpha$ is representable by (a multiple of) a closed geodesic. Moreover, if $K(p)=K<0$ is constant and
$\Sigma\approx\R^2$ then $\hat\kappa(m)$, $m\in\Z$, satisfies
\begin{equation*}
\hat\kappa(m)=
\begin{cases}
2\pi &\text{for $m=0$,}\\
\pi(|m|+1) &\text{for $m\ne 0$.}
\end{cases}
\end{equation*}
\item[{\rm (b)}]
Let $\Sigma$ be a Riemann surface of constant curvature $K=0$ (the completeness assumption then implies $\Sigma\approx \R^2$, $\Sigma\approx S^1\times\R$, or $\Sigma\approx T^2$). Then the infimum $\hat\kappa((\xi,m))$ is attained at some curve in the regular homotopy class
$(\xi,m)\in\pi_1(\Sigma)\times\Z$ if and only if $\xi\ne\ast$ or $m\ne 0$, where $\ast$ denotes the homotopy class of the constant loop. Moreover, $\hat\kappa((\xi,m))$, $(\xi,m)\in\pi_1(\Sigma)\times\Z$, satisfies
\begin{equation*}
\hat\kappa((\xi,m))=
\begin{cases}
2\pi &\text{for $(\xi,m)=(\ast,0)$,}\\
2\pi(|m|) &\text{otherwise.}
\end{cases}
\end{equation*}
\item[{\rm (c)}]
Let $\Sigma$ be the $2$--sphere with any metric and let $\alpha$ be a regular homotopy class of circle immersions into $\Sigma$. Then the infimum
$\hat\kappa(\alpha)$ equals $0$ and is attained at some curve in $\alpha$.
\end{itemize}
\end{thm}
\fullref{thminf} is proved in \fullref{ssec:5A}.

A curve in a Riemann surface $\Sigma$ with $K(p)\ne 0$ for all
$p\in\Sigma$, which is a local minimum of $\kappa$ is in fact a closed
geodesic, see \fullref{prploc->glob} (a). For flat Riemann
surfaces this is not the case. Here any local minimum of $\kappa$ is a
locally convex curve, see \fullref{prploc->glob} (b). We say
that a curve $c$ is {\em locally convex} if $k_g\ge 0$ everywhere for
some orientation of $c$. If $k_g>0$ everywhere, we say that $c$ is
{\em strictly locally convex}. In the terminology of Gromov \cite[page
8]{Gr}, strictly locally convex curves are called free curves.

The following result describes the homotopy types of the spaces of local minima of $\kappa$ for a flat Riemann surface. (Here we think of circle immersions as oriented unit speed curves parameterized by
arc length.)

\begin{thm}\label{thmfree}
On a flat Riemann surface $\Sigma$ (the completeness assumption implies $\Sigma\approx \R^2$, $\Sigma\approx S^1\times\R$, or $\Sigma\approx T^2$), the space $\Omega_{(\xi,m)}$of
(strictly) locally convex curves of regular homotopy class
$(\xi,m)\in\pi_1(\Sigma)\times\Z$ satisfies
$$
\Omega_{(\xi,m)}\simeq
\begin{cases}
\emptyset &\text{ if $(\xi,m)=(\ast,0)$,}\\
\Sigma &\text{if $\xi\ne\ast$ and $m=0$,}\\
U\Sigma &\text{if $m\ne 0$,}
\end{cases}
$$
where $\simeq$ denotes weak homotopy equivalence.
\end{thm}
\fullref{thmfree} is proved in \fullref{ssec:5B}.

\subsection*{Acknowledgements}
During the preparation of this paper, the author was a research fellow
of the Royal Swedish Academy of Sciences sponsored by the Knut and
Alice Wallenberg foundation. He also acknowledges support from the
Alfred P Sloan Foundation and from NSF--grant DMS-0505076. The author
thanks B Pausader for reading and commenting on earlier versions of
the paper.

\section{First variation and local minima}\label{sec:2}
In this section we compute the first variation of total absolute geodesic curvature. We use the result to classify local minima. Since the absolute value function is not differentiable at zero, the first variation is expressed as a statement about differences rather than as a statement about derivatives.

\subsection{First variation of total absolute geodesic curvature}\label{ssec:2A}
Let $c\co [0,L]\to \Sigma$ be an immersion into a Riemann surface, parameterized by arc length. Let $e_1$ be the unit tangent vector field of $c$ and let $e_2$ be a unit vector field along $c$ everywhere orthogonal to $e_1$. We consider variations $\omega\co[0,L]\times(-\delta,\delta)\to \Sigma$ of $c$ with the following three properties: $\omega(s,0)=c(s)$, $\pa_\epsilon \omega(s,0)=\alpha(s)e_2(s)$, where $\pa_\epsilon$ denotes differentiation with respect to the second variable, for some function $\alpha\co[0,L]\to\R$, and the curves $\omega_\epsilon(s)=\omega(s,\epsilon)$ ($\epsilon$ fixed) are immersions for $\epsilon\in(-\delta,\delta)$. We also introduce the sign function $\sigma\co\R\to\R$ as follows
\begin{equation*}
\sigma(x)=
\begin{cases}
1  &\text{ for }x>0,\\
0  &\text{ for }x=0,\\
-1 &\text{ for }x<0.
\end{cases}
\end{equation*}

\begin{lma}\label{lma1stvar}
For $\epsilon\in (-\delta,\delta)$,
\begin{align*}
\kappa(\omega_\epsilon)-\kappa(c) = & \epsilon\left(\int_{\{k_g\ne 0\}} \sigma(k_g)
\left(\ddot\alpha+\alpha K\right)\,ds\right)\\
& +|\epsilon|\left(\int_{\{k_g=0\}}|\ddot\alpha+\alpha K|\,ds\right) +
\Ordo(\epsilon^2),
\end{align*}
where $\ddot\alpha=\frac{d^2\alpha}{ds^2}$, and where $\Ordo(\epsilon^2)$ denotes a function such that $\epsilon^a\Ordo(\epsilon^2)\to 0$ as $\epsilon\to 0$ for all $a>-2$.
\end{lma}

\begin{proof}
To simplify notation, let $\frac{\pa \omega}{\pa s}=\dot \omega$ and $\frac{\pa \omega}{\pa \epsilon}=\omega'$. Let $\nabla$ denote the Levi--Civita connection and let $\nabla_s=\nabla_{\dot \omega}$ and $\nabla_\epsilon=\nabla_{\omega'}$. If $d\tau$ denotes the arc length element of the curve $\omega_\epsilon$ and $k_g(s,\epsilon)$ denotes the geodesic curvature of $\omega_\epsilon$ at $s$ then $$
|k_g|\,d\tau=
\left|k_g(s,\epsilon)\right|
\left|\dot \omega\right|\,ds=
\frac{\Bigl|\left\la\nabla_s\,\dot\omega,\iota\,\dot\omega\right\ra\Bigr|}
{\left|\dot \omega\right|^2}\,ds,
$$
where $\iota$ denotes rotation by $\frac{\pi}{2}$. Assuming $k_g(s,0)\ne 0$ and remembering that $c$ is parameterized by arc length, we compute
\begin{align*}
\pa_\epsilon\bigl(|k_g(s,0)||\dot \omega|\bigr)=&
\Bigl(\pa_\epsilon|\dot \omega|^{-2}\Bigr)|k_g|
+\sigma(k_g)
\Bigl(\pa_\epsilon\left\la\nabla_s\,\dot\omega,\iota\,\dot\omega\right\ra\Bigr)\\
=&-2\la\nabla_\epsilon\,\dot \omega,\dot \omega\ra|k_g|
+\sigma(k_g)
\Bigl(\left\la\nabla_\epsilon\nabla_s\,\dot\omega,\iota\,\dot\omega\right\ra
+\left\la\nabla_s\,\dot\omega,\iota\nabla_\epsilon\,\dot\omega\right\ra\Bigr)\\
=&-2\la\nabla_s\,\omega',\dot \omega\ra|k_g|\\
&+\sigma(k_g)\Bigl(\left\la\nabla_s\nabla_s\,\omega',\iota\,\dot\omega\right\ra
+\left\la R(\omega',\dot\omega)\dot\omega,\iota\,\dot\omega\right\ra
+\left\la\nabla_s\,\dot\omega,\iota\nabla_s\,\omega'\right\ra\Bigr),
\end{align*}
where $R$ is the curvature tensor. Noting that $\omega'=\alpha e_2$, $\nabla_s\, e_1=k_ge_2$, and $\nabla_s\, e_2=-k_ge_1$, we conclude
\begin{align}\notag
\pa_\epsilon|k_g(s,0)|&=
2\alpha k_g|k_g|+ \sigma(k_g)\Bigl(\ddot\alpha -\alpha k_g^2+ K\alpha
-\alpha k_g^2\Bigr)\\\label{eqder}
&= \sigma(k_g)\Bigl(\ddot\alpha + K\alpha\Bigr).
\end{align}
A similar calculation at $s$ where $k_g(s,0)=0$ gives
$$
\pa_\epsilon \bigl(k_g(s,0)|\dot \omega|\bigr)=\ddot\alpha + K\alpha.
$$
Hence for such $s$,
\begin{equation}\label{eqder1}
|k_g(s,\epsilon)||\dot\omega|=|\epsilon||\ddot\alpha +
 K\alpha|+\Ordo(\epsilon^2).
\end{equation}
The result follows by integration of \eqref{eqder} and \eqref{eqder1}.
\end{proof}

\subsection{Local minima of $\kappa$}\label{ssec:2B}
\begin{prp}\label{prploc->glob}
Let $\Sigma$ be a Riemann surface.
\begin{itemize}
\item[{\rm (a)}] If $K(p)\ne 0$ for all $p\in\Sigma$, then an
immersion $c\co S^1\to\Sigma$ is a local
minimum of $\kappa$ if and only if it is a geodesic.
\item[{\rm (b)}] If $\Sigma$ is flat then an
immersion $c\co S^1\to\Sigma$ is a local minimum of $\kappa$ if and only if it is
a locally convex curve.
\end{itemize}
In particular, in both cases {\rm (a)} and {\rm (b)}, any local
minimum of $\kappa$ is a global minimum in its regular homotopy class.
\end{prp}

\begin{proof}
Consider case (a). A curve $c\co S^1\to\Sigma$ is a geodesic if and only if $k_g(s)=0$ for all $s\in S^1$ and geodesics are global minima of $\kappa$. Let $c\co S^1\to\Sigma$ be a local minimum of $\kappa$.  Note that $U=\{s\in S^1\co k_g(s)\ne 0\}$ is open. Assume $U$ is nonempty. Then there exists a nonempty open subinterval $J\subset U$. For any variation $\omega(s,\epsilon)$ of $c$ with
$\pa_\epsilon \omega=\alpha e_2$ where $\alpha\co S^1\to\R$ is supported in $J$ we have
$$
\int_J\sigma(k_g)\ddot\alpha\,ds=0.
$$
Thus, since $c$ is a local minimum, we conclude from \fullref{lma1stvar} that
$$
\int_J\sigma(k_g)\alpha K\,ds=0.
$$
This contradicts $K(p)\ne 0$ for all $p\in\Sigma$. It follows that
$U$ is empty and thus $c$ is a geodesic.

Consider case (b). Let $c\co S^1\to\Sigma$ be a local minimum. We show that $k_g$ cannot change sign along $c$. Assume it does, then there exist two disjoint open subintervals $J_+$ and $J_-$ of $S^1$ such that $k_g>0$ on $J_+$ and $k_g<0$ on $J_-$. Let $A$ be a subinterval of $S^1$ containing both $J_+$ and $J_-$. Let $\alpha\co A\to\R$ be a function such that $\ddot\alpha$ is supported in small subintervals of $J_+\cup J_-$ and such that $\dot\alpha=r\ne 0$, where $r$ is a non-zero constant, between $J_+$ and $J_-$. For a variation $\omega$ of $c$ with $\pa_\epsilon\omega=\alpha e_2$, Lemma 2.1 implies that
$$
\kappa(\omega_\epsilon)-\kappa(c)=\pm 2r\epsilon+\Ordo(\epsilon^2).
$$
This contradicts $c$ being a local minimum. Consequently, $k_g$ does not change sign along $c$ and $c$ is locally convex.

It remains to show that $c$ is a global minimum. Fix a unit speed parametrization of $c$ so that $k_g(s)\ge 0$ for all $s\in S^1$. As in \fullref{sec:1}, we construct an orthonormal trivialization of $T\Sigma$ by parallel translation with respect to the flat metric. This identifies the unit tangent bundle $U\Sigma$ of $\Sigma$ with $\Sigma\times S^1$ and the regular homotopy class of $c$ is determined by its homotopy class in $\Sigma$ and the degree of $\pi_2\circ\vec c\co S^1\to S^1$, where $\pi_2\co \Sigma\times S^1\to S^1$. (Recall that $\vec c\co S^1\to U\Sigma$ denotes the natural lift of the unit speed curve $c$.) Moreover, $\kappa(c)$ is simply the length of the curve $\pi_2\circ \vec c$. Now, $k_g(s)\ge 0$ for all $s\in S^1$ implies that the length of $\pi_2\circ \vec  c$ equals $2\pi$ times the degree of $\pi_2\circ \vec c$ and it follows that local minima are global minima also in this case.
\end{proof}

\begin{rmk}
\fullref{prploc->glob} does not hold for arbitrary Riemann surfaces. Consider for example the boundary of a convex body in $\R^3$ which agrees with the standard $2$--sphere except that it has a flat region near the north pole. Any locally convex curve in this flat region is a local minimum of $\kappa$ but it is certainly not a global minimum in its regular homotopy class.
\end{rmk}

\section{Curvature concentrations and approximations}\label{sec:3}
In this section we define piecewise geodesic curves with curvature concentrations and show that circle immersions can be approximated by such curves without increasing the total absolute geodesic curvature.

\subsection{Piecewise geodesic curves with curvature concentrations}\label{ssec:3A}
Let $\Sigma$ be a Riemann surface. A {\em piecewise geodesic curve} in $\Sigma$ is a continuous curve $c\co S^1\to\Sigma$ which is a finite union of geodesic segments. More formally such a curve $c$ can be described as follows. Consider a finite collection of geodesics $c_j\co [0,1]\to\Sigma$, $j=1,\dots,m$, with $c_j(1)=c_{j+1}(0)$ for each $j$ (here $c_{m+1}=c_1$). Let $I_j$ be $[0,1]$ thought of as the domain of $c_j$ and let $I_{m+1}=I_1$. Then the space obtained by identifying $1\in I_j$ with $0\in I_{j+1}$ is a circle $S^1$ which can be considered as the domain of a continuous map $c\co S^1\to\Sigma$ such that if $p\in S^1$ is the image of $p'\in I_j$ under the quotient projection then $c(p)=c_j(p')$. We say that the points $p\in S^1$ with two preimages under the quotient projection are the {\em vertices} of the piecewise geodesic curve $c$. We will often deal with images of vertices and we call also these image points vertices of $c$, when no confusion can arise.

Let $U\Sigma$ denote the unit tangent bundle of $\Sigma$. Note that at each vertex $c(p)$ of a piecewise geodesic curve $c$ as above, there is an {\em incoming} unit tangent vector ${\dot c}_j(1)/|{\dot c}_j(1)|\in U_{c(p)}\Sigma$ and an {\em outgoing} unit tangent vector ${\dot c}_{j+1}(0)/|{\dot c}_{j+1}(0)|\in U_{c(p)}\Sigma$. A piecewise geodesic curve {\em with curvature concentrations} is a piecewise geodesic curve $c$ together with a {\em vertex curve} $\gamma\co[0,1]\to U_{c(p)}\Sigma$ for each vertex $p$ which connects the incoming-- to the outgoing unit tangent of $c$ at $p$ and which satisfies the following condition: $\gamma$ is a continuous piecewise geodesic curve (with finitely many geodesic arcs) in the fiber circle $U_{c(p)}\Sigma$ equipped with the metric induced by the Riemannian metric on $\Sigma$. We use the abbreviations $\pgc$--curve to denote piecewise geodesic curves with curvature concentrations, and we often write $(c_1,\gamma_1,\dots,c_m,\gamma_m)$ for a $\pgc$--curve with geodesic segments $c_j$ and vertex curves $\gamma_j$. We say that the length $l(\gamma_j)$ of the vertex curve $\gamma_j$ is the {\em curvature concentration} of the $\pgc$--curve $c$ at the vertex $c(p)=c_j(1)=c_{j+1}(0)$, and that the piecewise geodesic curve with geodesic arcs $c_1,\dots,c_m$ is the {\em underlying curve} of $c$.

We note that any $\pgc$--curve $c=(c_1,\gamma_1,\dots, c_m,\gamma_m)$ in $\Sigma$ has a natural continuous lift $\vec c\co S^1\to U\Sigma$ which consists of the usual lifts $\vec c_j$ of $c_j$, $j=1,\dots,m$, connected by the curves $\gamma_j$, $j=1,\dots,m$, in the fibers of $U\Sigma\to\Sigma$ over vertices of $c$. The lift $\vec c$ of $c$ is thus a piecewise smooth curve. In particular, its derivative is smooth except for finitely many jump discontinuities where the curve has left and right derivatives. If $c$ is a $\pgc$--curve we consider $\vec c\co  S^1\to U\Sigma$ as a parameterized curve with its natural arc length parametrization scaled by a suitable factor so that its domain becomes the unit circle.

\begin{dfn}
The {\em total absolute geodesic curvature} of a $\pgc$--curve
$$
c=(c_1,\gamma_1,\dots,c_m,\gamma_m)
$$
$$
\kappa(c)=\sum_{j=1}^m l(\gamma_j),\leqno{\hbox{is}}
$$
where $l(\gamma_j)$ is the length of the vertex curve $\gamma_j$.
\end{dfn}

Let $\pImm(S^1,\Sigma)$ denote the set of all $\pgc$--curves in $\Sigma$. We define the distance between two elements $b$ and $c$ in $\pImm(S^1,\Sigma)$ to be the $C^0$--distance (with respect to the metric on $U\Sigma$ induced from the metric on $\Sigma$) between their lifts $\vec b$ and $\vec c$ endowed with parameterizations proportional to arc length, as discussed above. A {\em $\pgc$--homotopy} is a continuous $1$--parameter family of $\pgc$--curves or equivalently a continuous map from the interval to $\pImm(S^1,\Sigma)$.

\subsection{Approximation}\label{ssec:3B}
If $\Sigma$ is a Riemann surface then let $\Imm(S^1,\Sigma)$ denote the space of circle immersions into $\Sigma$ with the $C^2$--topology. Let $c\co S^1\to\Sigma$ be a circle immersion, let $\pi=(p_0,\dots,p_m)$ be a partition of $S^1$ and let $|\pi|=\max_jd(p_j,p_{j+1})$ ($d$ is the distance function on $S^1$ and we use the convention $p_{m+1}=p_0$). If $|\pi|$ is sufficiently small then we associate a $\pgc$--curve $c^\pi$ to $c$, as follows: $c^\pi$ is the $\pgc$--curve with underlying piecewise geodesic curve consisting of the shortest geodesic segments between $c(p_j)$ and $c(p_{j+1})$, and with vertex curves $\gamma_j$ which are the shortest arcs in $U_{c(p_j)}\Sigma$ connecting the incoming-- to the outgoing unit tangent of the underlying piecewise geodesic curve. We note that $\max_j l(\gamma_j)\to 0$ and $|\kappa(c)-\kappa(c^\pi)|\to 0$ as $|\pi|\to 0$.

If $f\co\Lambda\to\Imm(S^1,\Sigma)$ is a continuous family of circle immersions parameterized by a compact space $\Lambda$ and if $\epsilon>0$ is arbitrary, then there exists $\delta>0$ such that for all partitions $\pi$ with $|\pi|<\delta$, $f^\pi\co\Lambda\to\pImm(S^1,\Sigma)$, defined by $f^\pi(\lambda)=(f(\lambda))^\pi$, is a continuous family of $\pgc$--curves, and there exists an $\epsilon$--small homotopy connecting the family of continuous curves $\vec f\co\Lambda\to U\Sigma$ to the family of continuous curves $\vec {f^\pi}\co\Lambda\to U\Sigma$.

\begin{lma}\label{lmasmooth->pgc}
Let $c_0\co S^1\to\Sigma$ be a circle immersion into a Riemann surface with $K<0$, $K>0$, or $K=0$ everywhere. Then there exists a regular homotopy $c_t$, $0\le t\le 1$, of $c_0$ such that
\begin{equation}\label{eqapp}
\kappa(c_t)\le\kappa(c_0),\quad 0<t\le 1,
\end{equation}
and such that the $\pgc$--curve $c_1^\pi$ (defined using any sufficiently fine partition $\pi$), satisfies
\begin{equation}\label{eqapp'}
\kappa(c_1^\pi)\le\kappa(c_0).
\end{equation}
Moreover, if $c$ is not a local minimum of $\kappa$ then the non-strict inequalities in \eqref{eqapp} and \eqref{eqapp'} can be replaced by strict inequalities.
\end{lma}

\begin{proof}
Assume $K\ne 0$ everywhere. Then \fullref{prploc->glob} implies that $c_0$ is a local minimum if and
only if $c_0$ is a geodesic. Hence, if $c_0$ is not a geodesic then there exists a $\kappa$--decreasing regular homotopy connecting $c_0$ to some curve $c_1$. For sufficiently fine partition $\pi$, the  curve $c_1^\pi$ then satisfies \eqref{eqapp'}.

Assume $K=0$. If $c$ is not locally convex then the above argument can be repeated. Recall from \fullref{sec:1} that parallel translation in the flat metric gives $U\Sigma=\Sigma\times S^1$ and let $\pi_2\co U\Sigma\to S^1$ denote the projection. If $c$ is locally convex it is elementary to see that for sufficiently fine partitions $\pi$, $\pi_2\circ \vec{c^\pi}$ is monotone and thus $\kappa(c^\pi)=\kappa(c)$.
\end{proof}

\begin{rmk}
\fullref{lmasmooth->pgc} does not hold for general metrics. Consider for example $\R^2$ with coordinates $(x,y)$ and a metric given by
$$
ds^2=\exp\Bigl(2\left(\sqrt{x^2+y^2}-1\right)^3\Bigr)(dx^2+dy^2)
$$
in a small neighborhood of $c=\{x^2+y^2=1\}$. Then $K(x,y)<0$ for $x^2+y^2>1$ and $K(x,y)>0$ for $x^2+y^2<1$ and the geodesic
curvature of $c$ is identically equal to $1$. Let $b$ be any curve which is a $C^1$--small perturbation of $c$. Assume that $b$ meets $c$ transversely in $2m$ points. These intersection points subdivide $b$ and $c$ into unions of arcs $c=\gamma_{{\rm out}}\cup\gamma_{{\rm in}}$ and $b=\beta_{{\rm out}}\cup\beta_{{\rm in}}$, where $\beta_{\rm out}\subset\{x^2+y^2\ge 1\}$ and $\beta_{{\rm in}}\subset\{x^2+y^2\le 1\}$ and where the endpoints of an arc in $\gamma_{{\rm in}}$ ($\gamma_{{\rm out}}$) agree with the endpoints of some arc in $\beta_{{\rm out}}$ ($\beta_{{\rm in}}$). Let $\Omega_{{\rm out}}^j$, $j=1,\dots, m$, be the $m$ regions bounded by an arc in $\beta_{{\rm out}}$ and an arc in $\gamma_{{\rm in}}$ and let $\Omega_{{\rm in}}^j$, $j=1,\dots,m$ be the $m$ region bounded by an arc in $\gamma_{{\rm out}}$ and an arc in $\beta_{{\rm in}}$. Then the Gauss--Bonnet theorem implies that
\begin{align*}
&\int_{\beta_{{\rm out}}}k_g\, ds - \int_{\gamma_{\rm in}}k_g\, ds +\sum_{j=1}^m\int_{\Omega_{{\rm out}}^j} K\,dA+\sum_{j=1}^{2m}\alpha_j=2\pi m,\\
&\int_{\gamma_{{\rm out}}}k_g\, ds - \int_{\beta_{\rm in}}k_g\, ds +\sum_{j=1}^m\int_{\Omega_{{\rm in}}^j} K\,dA+\sum_{j=1}^{2m}\alpha_j=2\pi m,
\end{align*}
where $\alpha_j$ is the exterior angle at the $j^{\rm th}$ intersection point between $b$ and $c$. Thus
$$
\kappa(b)\ge \int_b k_g\, ds > \int_c k_g\, ds=\kappa(c).
$$
By approximation we conclude that for any $\pgc$--curve $e$ which lies in a sufficiently small tubular neighborhood of $c$, $\kappa(e)>\kappa(c)$.
\end{rmk}

\section{Curvature non-increasing homotopies, smoothing, and locally convex curves}\label{sec:4}
In this section we construct special $\pgc$--homotopies of $\pgc$--curves on flat Riemann surfaces and on the hyperbolic plane which decrease the total absolute geodesic curvature of a given initial curve and which ends at a curve of certain standard shape. We construct similar special $\pgc$--homotopies of curves on the $2$--sphere with a constant curvature metric. We show that these special $\pgc$--homotopies can be smoothed to regular homotopies, increasing the total curvature arbitrarily little. We also study the space of locally convex curves on flat surfaces.

\subsection{Flat surfaces}\label{ssec:flat}\label{ssec:4A}
Let $\Sigma$ be a Riemann surface. Let $\Pi\co\tilde\Sigma\to\Sigma$ be a smooth covering map and endow $\tilde\Sigma$ with the pull-back metric. Define a {\em lift} of a $\pgc$--curve $c$ in $\Sigma$ to be a $\pgc$--curve $b$ in $\tilde\Sigma$ such that $\vec b$ is a lift of $\vec c$ with respect to the induced covering $\Pi_U\co U\tilde\Sigma\to U\Sigma$.

Let $c$ be a $\pgc$--curve in $\Sigma$ with lift $b$ in $\tilde\Sigma$ and let $b_t$, $0\le t\le 1$, be a $\pgc$--homotopy of $b$ with the following properties: the start-point $b_t(\alpha)$ and endpoint $b_t(\omega)$ of $b_t$ satisfy $\Pi(b_t(\alpha))=\Pi(b_t(\omega))$ and $\Pi_U(\dot b_t(\alpha))=\Pi_U(\dot b_t(\omega))$, for all $t$, where $\dot b_t(\alpha)$ and $\dot b_t(\omega)$ are the outgoing-- and incoming tangent vectors of $b_t$ at $b_t(\alpha)$ and $b_t(\omega)$, respectively. Then $b_t$ induces a $\pgc$--homotopy $c_t$ of $c$ by transporting geodesic segments of $b_t$ in $\tilde\Sigma$ to geodesic segments of $c_t$ in $\Sigma$ with the projection $\Pi\co \tilde\Sigma\to\Sigma$, and by transporting vertex curves of $b_t$ to vertex curves of $c_t$ with the induced projection $\Pi_U\co U\tilde\Sigma\to U\Sigma$. (Note that the conditions on the $\pgc$--homotopy $b_t$ holds if $\vec b_t$ is a closed curve for each $t$.)

\begin{lma}\label{lmaflatV}
Let $\Sigma$ be a Riemann surface with constant curvature $K=0$ and let $c$ be a $\pgc$--curve in $\Sigma$. Then there exists a $\pgc$--homotopy $c_t$, $0\le t\le 1$, with $c_0=c$, with $\kappa(c_t)\le\kappa(c_0)$, for all $t$, and with the following property: the underlying curve of $c_1$ is a geodesic and any curvature concentration of $c_1$ equals $\pi$.
\end{lma}

\begin{proof}
Consider the case $\Sigma=\R^2$. We claim that any $\pgc$--curve which does not have underlying curve a line segment admits a $\pgc$--homotopy which does not increase $\kappa$ and which decreases the number of vertices with curvature concentration not a multiple of $\pi$. Together with an obvious inductive argument this shows that any $\pgc$--curve is $\pgc$--homotopic through an homotopy with properties as above to a curve with all curvature concentrations integral multiples of $\pi$. Noting that any curvature concentration at a vertex of magnitude $m\pi$, $m$ a positive integer, can be split up into $m$ vertices each with curvature concentration $\pi$ by a $\pgc$--homotopy preserving $\kappa$ we conclude that the claim implies the lemma when $\Sigma=\R^2$.

Consider the claim. If the underlying curve of a $\pgc$--curve is not a line segment then it has three consecutive line segments connected by two vertices with curvature concentrations which are not integral multiples of $\pi$. At such a vertex the curvature concentration is either the sum of the exterior angle of the underlying curve and a multiple of $2\pi$ or the sum of the interior angle of the curve and an odd multiple of $\pi$. We call the former type of vertices {\em exterior} and the latter {\em interior}. We also need to distinguish two types of configurations of the three consecutive segments. We say that the configuration is {\em convex} if all three segments are contained in the closure of one of the half planes determined by the line containing the middle segment, otherwise we say it is {\em non-convex}. To establish the claim, consider three consecutive segments $e_1$, $e_2$, $e_3$ as above connected at vertices $v_1$ and $v_2$. We separate the cases:

{\bf Case 1}\qua If both $v_1$ and $v_2$ are exterior, then move $v_2$ along $e_3$ until it reaches the next vertex following it. It is straightforward to check that this $\pgc$--homotopy preserves $\kappa$ in the convex case and decreases $\kappa$ in the non-convex case.

{\bf Case 2}\qua If $v_1$ is exterior and $v_2$ interior, then move $v_1$ backwards along $e_1$ until it reaches the vertex preceding it. This $\pgc$--homotopy preserves $\kappa$ in the non-convex case and decreases it in the convex case.

{\bf Case 3}\qua Assume that both $v_1$ and $v_2$ are interior. If the configuration is convex, then move $v_2$ along $e_3$ until it reaches the vertex following it. This is a $\kappa$--preserving $\pgc$--homotopy. If the configuration is non-convex, consider the lines $l_j$ containing $e_j$. Assume first that $l_1$ and $l_3$ are not parallel. Note that the segments $e_1$ and $e_3$ lie in different components of $\R^2-l_2$. If the point $l_1\cap l_3$ lies in the component of $\R^2-l_2$ which contains  $e_3$ then move $v_1$ along $l_1$ to $l_1\cap l_3$. Otherwise, move $v_2$ along $l_3$ to the intersection point. This $\pgc$--homotopy strictly decreases $\kappa$. In the case that $l_1$ and $l_3$ are parallel, start by moving $v_1$ as described above. Note that this decreases $\kappa$. Thus we may change $e_3$ slightly without increasing $\kappa$ past $\kappa(c_0)$ so that $l_1$ and $l_3$ intersect. We then apply the above.

Note that in either case, the $\pgc$--homotopy described reduces the number of vertices with curvature concentration not an integral multiple of $\pi$. The claim follows.

Consider the non-simply connected case. Since $\Sigma$ is a flat Riemann surface, there is a covering map $\Pi\co\R^2\to\Sigma$ with deck transformations which are translations of $\R^2$. Let $b$ be a lift of $c$. If $b$ is closed, we apply the above result to construct a $\pgc$--homotopy $b_t$, $0\le t\le 1$, with properties as in the formulation of the lemma. The induced $\pgc$--homotopy $c_t$, $0\le t\le 1$, of $c$ then satisfies the lemma. If $b$ is non-closed we lift it in the middle of some segment and apply the construction above to the lift. If at some stage of the $\pgc$--homotopy a moving vertex passes the endpoint of the lifted curve, we stop the $\pgc$--homotopy and choose the midpoint of the next segment of the curve in the direction the vertex is moving and then proceed. Note that this construction applies as long as the underlying curve has at least three distinct geodesic segments. Thus, to complete the proof it remains only to study curves with two line segments and two vertices.

Let $b$ be the lift of a curve $c$ with two vertices. Let $v_0$ and $v_1$ be the endpoints of the lift $b$ and let $l$ be the straight line segment in $\R^2$ connecting $v_0$ to $v_1$. Let $e_0$ and $e_1$ be the straight line segments of $b$ with one endpoint at $v_0$ and at $v_1$, respectively. Let $e$ be the segment of $b$ connecting the endpoint of $e_0$ which is not equal to $v_0$ to the endpoint of $e_1$ not equal to $v_1$. Since we lift $c$ at a midpoint of a geodesic segment it follows that $e$ passes the midpoint $p$ of $l$.

Let $w_0$ and $w_1$ be the vertices where $e_0$ and $e$ meet and where $e_1$ and $e$ meet, respectively. Choose an orienting basis $(\hat e,\hat f)$ of the plane where $\hat e$ is a unit vector in direction of $e$ oriented from $v_0$ to $v_1$ and where $\hat f$ is a vector perpendicular to $\hat e$ oriented so that the direction vector of $e_0$ has positive $\hat f$ component. If $\Delta(r,s,t)$ is a triangle with corners at $r,s,t\in\R^2$ we write $\alpha(r;r,s,t)$ and $\beta(r;r,s,t)$ for the exterior-- respectively interior angle of $\Delta(r,s,t)$ at the corner $r$.

We consider separate cases. Consider first the case when both $w_0$ and $w_1$ are exterior vertices. Rotate the line containing $e$ in the positive directions around $p$, and rotate the lines containing $e_0$ and $e_1$ in the negative direction around $v_0$ and $v_1$, respectively, so that angles change at linear speed. More precisely if $w_0^t$ and $w_1^t$ denotes the intersection points between the rotated line containing $e$ and the rotated line containing $e_0$ and $e_1$, respectively. Then $\beta(p;p,w_j^t,v_j)=(1-t)\beta(p;p,w_j,v_j)$, $j=0,1$. Note that the segments in the triangles $\Delta(p,w_j^t,v_j)$, $j=0,1$, which are not parallel to $l$ give a $\pgc$--homotopy $b_t$ of $b$ which via the projection $\Pi$ induces a $\pgc$--homotopy $c_t$ of $c$. Moreover,
$$
\frac{d}{dt}\Bigl(\kappa(c)-\kappa(c_t)\Bigr)=2\bigl(\beta(p;p,v_0,w_0)+\beta(v_0;p,v_0,w_0)\bigr)>0,
$$
and the curvature concentrations of $c_1$ are integral multiples of $2\pi$.

Consider second the case when both $w_0$ and $w_1$ are interior angles. In this case we repeat the construction above with the important difference that we rotate also the lines containing $e_0$ and $e_1$ in the positive direction.  Again the projection gives a $\pgc$--homotopy $c_t$ of $c$ and we have
$$
\frac{d}{dt}\Bigl(\kappa(c_t)-\kappa(c)\Bigr)=
2\bigl(-\alpha(v_0;p,v_0,w_0)+\beta(p;p,v_0,w_0)\bigr),
$$
but $\alpha(v_0;p,v_0,w_0)=\beta(p;p,v_0,w_0)+\beta(w_0;p,v_0,w_0)$ and hence $\kappa(c_t)<\kappa(c)$. Again $c_1$ has the desired form.

Finally, if one of $w_0$ and $w_1$ is an exterior vertex and the other one is an interior vertex, then any one of the above procedures may be used. The result is a $\pgc$--homotopy which does not change $\kappa$: the sum of the interior angle and the exterior angle at $w_0$ and $w_1$ is constantly equal to $\pi$ . This finishes the proof in the non-simply connected case.
\end{proof}

\begin{lma}\label{lmaflat>}
If $c_0$ and $c_1$ are two $\pgc$--homotopic $\pgc$--curves on a Riemann surface as in \fullref{lmaflatV} then there exists a $\pgc$--homotopy $c_t$, $0\le t\le 1$, connecting them such that
$$
\kappa(c_t)\le\max\{\kappa(c_0),\kappa(c_1)\},
$$
for all $t$.
\end{lma}

\begin{proof}
After \fullref{lmaflatV} it is sufficient to consider the case when $c_0$ and $c_1$ both have underlying curves geodesics and all curvature concentrations equal to $\pi$. Any vertex curve is thus either a positive or a negative $\pi$--rotation. It is easy to see that neighboring vertex curves of different orientations cancel. The lemma follows.
\end{proof}

\subsection{The hyperbolic plane}\label{ssec:4B}
The counterpart of \fullref{lmaflatV} for curves in the hyperbolic plane differs from the flat case in an essential way: the limit curve which arises as the end of a $\kappa$--decreasing $\pgc$--homotopy, is not a $\pgc$--curve {\em in} the hyperbolic plane, in fact it often has infinite length. To deal with this phenomenon we define a generalized $\pgc$--curve in the hyperbolic plane as a $\pgc$--curve which is allowed to have vertices at infinity. More concretely, consider the disk model of the hyperbolic plane,
\begin{equation}\label{eqH2}
D=\{x=(x_1,x_2)\in\R^2\co |x|<1\},\quad
ds^2=\frac{4(dx_1^2+dx_2^2)}{(1-|x|^2)^2},
\end{equation}
and add to it the circle at infinity $\pa D$. Define a {\em generalized $\pgc$--curve} as a piecewise smooth curve in $\bar D=D\cup\pa D$ which consists of geodesic segments (ie, arcs of circles meeting $\pa D$ at right angles), which is allowed to have vertices on $\pa D$, and which have vertex curves at all vertices connecting the incoming-- to the outgoing unit tangent. Note that the length of any vertex curve at a vertex on $\pa D$ is an odd multiple of $\pi$. We extend $\kappa$ to generalized $\pgc$--curves by defining it as the sum of the lengths of all vertex curves (the sum of all curvature concentrations).

To connect generalized $\pgc$--curves to $\pgc$--curves we measure the distance between generalized $\pgc$--curves as the $C^0$--distance between their lifts in $U \bar D$ with respect to the metric on $U \bar D$ induced by the Euclidean metric on the plane. Using this metric we define a generalized $\pgc$--homotopy as a continuous $1$--parameter family of generalized $\pgc$--curves. Moreover, it is clear that the following approximation result holds: if $b_\lambda$, $\lambda\in\Lambda$ is any continuous family of generalized $\pgc$--curves parameterized by a compact space $\Lambda$ and if $\epsilon>0$ is arbitrary then there exists a family of $\pgc$--curves $c_\lambda$, $\lambda\in\Lambda$, such that the distance between $b_\lambda$ and $c_\lambda$ is less than $\epsilon$ and such that $|\kappa(b_\lambda)-\kappa(c_\lambda)|<\epsilon$ for every $\lambda\in\Lambda$.

\begin{lma}\label{lmanegV}
Let $\Sigma\approx\R^2$ be a Riemann surface of constant curvature $K<0$ and let $c$ be a $\pgc$--curve in $\Sigma$. Then there exists a generalized $\pgc$--homotopy $c_t$, $0\le t\le 1$, with $c_0=c$, with $\kappa(c_t)$ a non-increasing function of $t$, and with the following property: $c_1$ is a generalized $\pgc$--curve with all its vertices on the circle at infinity, and with all curvature concentrations equal to $\pi$.
\end{lma}

\begin{proof}
Note that $\kappa$ is invariant under scaling. Therefore we may assume that $K=-1$. Then $\Sigma=D$, where $D$ is as in \eqref{eqH2}. Let $c$ be a $\pgc$--curve in $\Sigma$. Consider a vertex curve $\gamma$ of $c$ at the vertex $p\in D$. Let $e_1$ be the incoming geodesic arc of $c$ at $p$. We consider two cases separately.

{\bf Case 1}\qua Assume that $l(\gamma)=a+n2\pi$ where $0< a <\pi$ and $n\ge 0$. In this case the vertex curve is, modulo $2\pi$--rotations, the exterior angle of the curve. We push the endpoint of $e_1$ backwards along $e_1$ until we reach its
start-point. At this moment we have reduced the number of non-infinite vertices of $c$ with curvature concentrations not a multiple of $\pi$ by one. To see that this generalized $\pgc$--homotopy does not increase $\kappa$ we calculate with notation as in \fullref{fig:neg1},
$$
\alpha'+\beta'=\alpha+\int_\Omega K\,dA\le\alpha,
$$
and hence $\kappa$ does not increase.

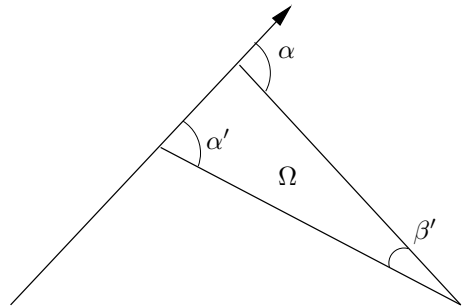
\begin{figure}[ht!]\centering
\begin{overpic}[angle=-90,width=6cm]{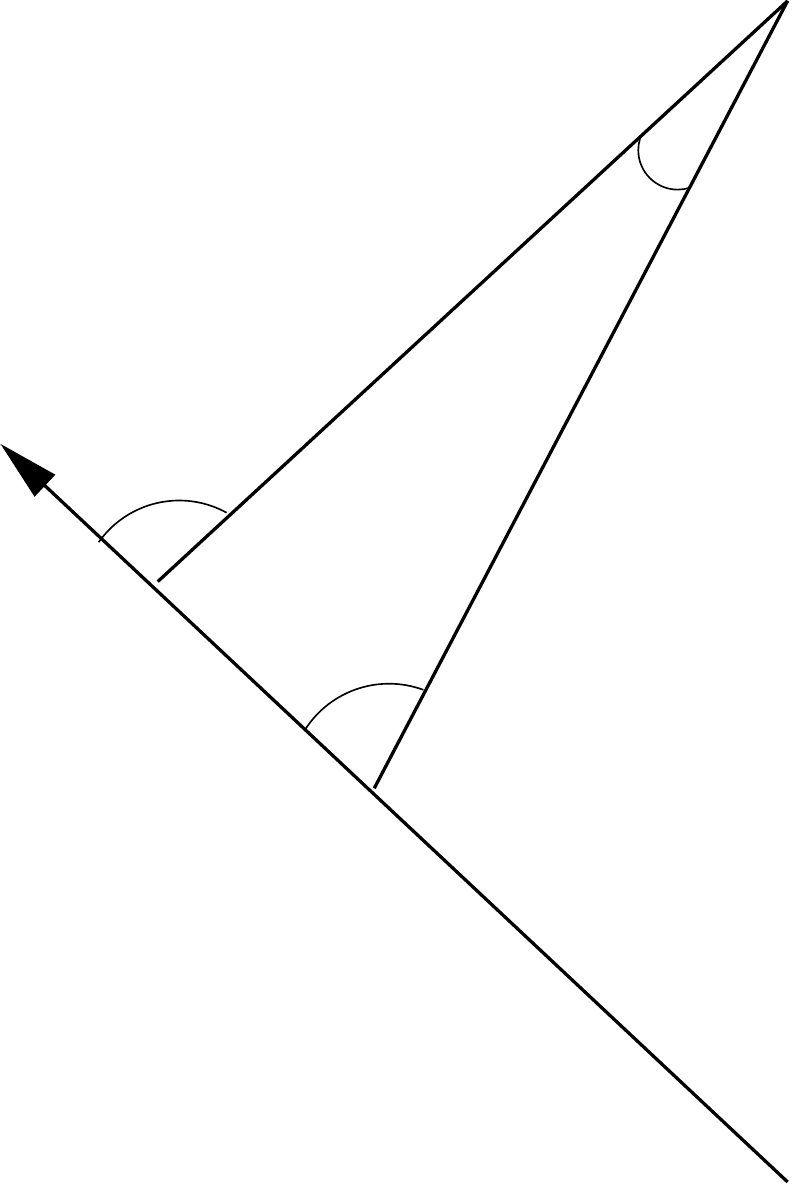}
{\small
\put(59,55){$\alpha$}
\put(43,35){$\alpha'$}
\put(59,27){$\Omega$}
\put(89,15){$\beta'$}
}
\end{overpic}
\caption{Removing an exterior angle}
\label{fig:neg1}
\end{figure}

{\bf Case 2}\qua Assume that $l(\gamma)=a+n2\pi$ where $\pi\le a<2\pi$ and $n\ge 0$. In this case the vertex curve is,
modulo $2\pi$--rotations, the complementary angle to the exterior angle. We push the vertex in the positive direction along $e_1$ until it hits $\pa D$, thereby reducing the number of non-infinite vertices by one. To see that this $\pgc$--homotopy does not increase $\kappa$ we calculate with notation as in \fullref{fig:neg2},
$$
\alpha'+\beta'=\alpha+\int_\Omega K\,dA\le\alpha,
$$
and hence $\kappa$ does not increase.

\begin{figure}[ht!]
\begin{overpic}[angle=-90,width=6cm]{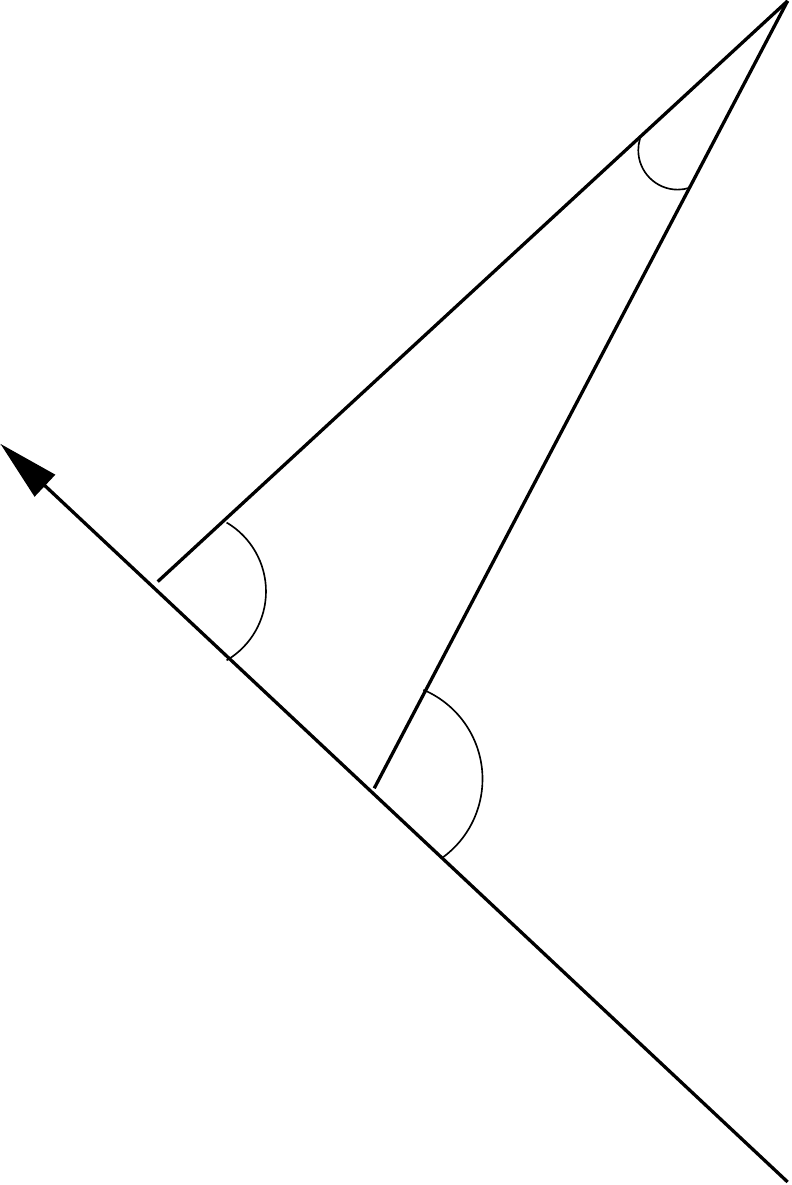}
{\small
\put(32,22){$\alpha$}
\put(47,39){$\alpha'$}
\put(59,27){$\Omega$}
\put(91,12){$\beta'$}
}
\end{overpic}
\caption{Removing an interior angle}
\label{fig:neg2}
\end{figure}
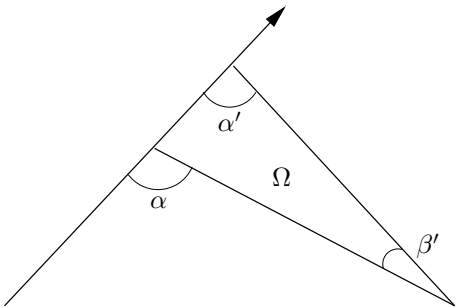

Repeating this argument a finite number of times we remove all non-infinite vertices of $c$. To finish the proof we note that any vertex curve at infinity has length $\pi + 2m\pi$, for some integer $m\ge 0$, and that the $2\pi m$--concentration can be pushed along one of the geodesics which ends at the infinite vertex so that it lies in $D$ (not in $\pa D$). Finally, such a curvature concentration of magnitude $2\pi m$ in the finite part of the disk can be split up and pushed to $2m$ curvature concentrations at infinity, each of length $\pi$. In  \fullref{fig:neg4} this is illustrated for $m=1$. (If the multiplicity of the geodesic on the left hand picture in \fullref{fig:neg4} is $k$, then the multiplicity between the two curvature concentrations on the right hand side is $k+2$.)
\begin{figure}[ht!]
\begin{center}
\begin{overpic}[angle=-90,width=8cm]{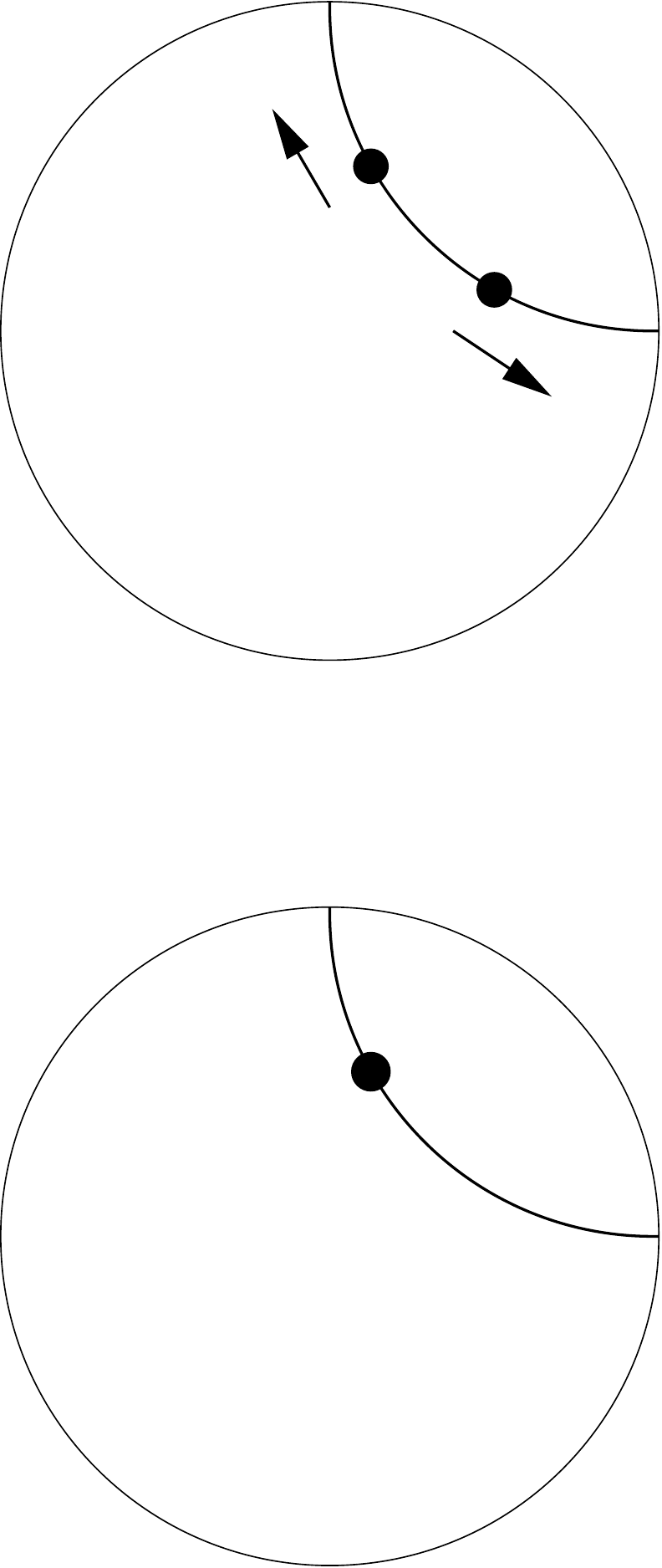}
{\small
\put(25,20){$2\pi$}
\put(83,8){$\pi$}
\put(90,15){$\pi$}
}
\end{overpic}
\end{center}
\caption{Splitting a curvature concentration}
\label{fig:neg4}
\end{figure}
\end{proof}

Our next goal is to deform any generalized $\pgc$--curve with its vertices at infinity to a standard form. To this end, we distinguish two different vertices at infinity: let $x$ be a vertex on $\pa D$ of a generalized $\pgc$--curve with incoming geodesic segment along the geodesic  $c_i$ and outgoing geodesic segment along the geodesic $c_o$. Note that $c_i$ subdivides $D$ into two components. Let $D_+$ be the component with inward normal $\nu$ along $c_i$ such that if $n(x)$ is the outward normal of $\pa D$ at $x$ then $n(x),\nu(x)$ is a positively oriented basis of $\R^2$. Let $D_-$ be the other component. If $c_i=c_o$ then we say $x$ is a {\em degenerate vertex}. Assume that $x$ is a non-degenerate vertex and that the vertex curve at $x$ is a positive (negative) $\pi$--rotation. Then we say that $x$ is an {\em over-rotated vertex} if $c_o$ lies in $D_-$ ($D_+$), otherwise we say it is an {\em under-rotated vertex}.

\begin{figure}[ht!]
\begin{center}
\includegraphics[angle=-90, width=5.5cm]{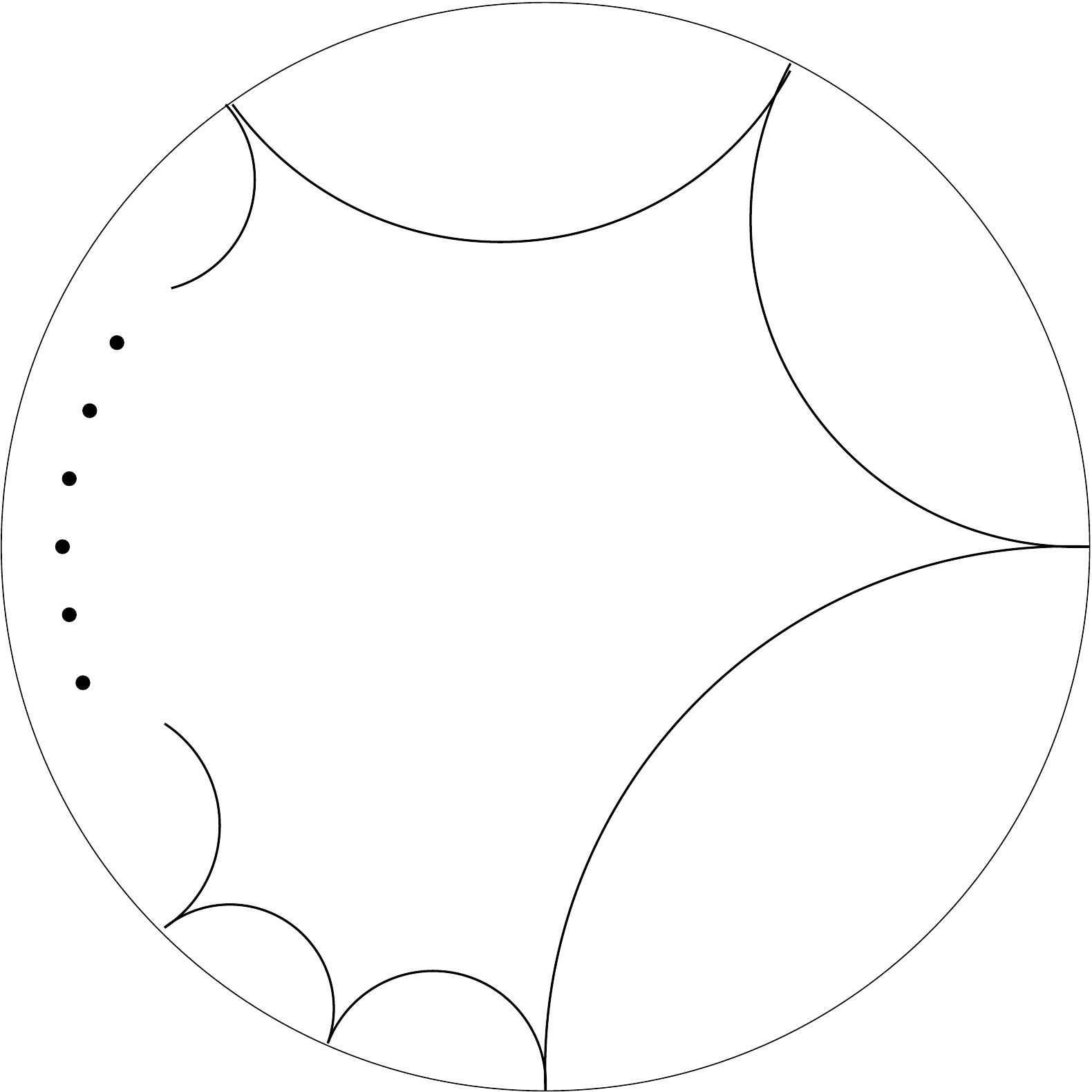}
\end{center}
\caption{A curve in standard position}
\label{fig:neg3}
\end{figure}

We say that a generalized $\pgc$--curve of tangential degree (Whitney index) $m$ in the hyperbolic plane is in {\em standard position} if it has the form of the curve in \fullref{fig:neg3}, where if, $|m|\ne 0$, all vertex curves have length $\pi$ and have the same orientation. Clearly any two $\pgc$--homotopic generalized $\pgc$--curves in standard position are $\pgc$--homotopic through such curves.

\begin{lma}\label{lmaneg>}
If $c_0$ and $c_1$ are two $\pgc$--homotopic $\pgc$--curves in the hyperbolic plane each with at least one curvature concentration not an integral multiple of $\pi$, then there exists a $\pgc$--homotopy $c_t$, $0\le t\le 1$, from $c_0$ to $c_1$ such that
$$
\kappa(c_t)\le\max\{\kappa(c_0),\kappa(c_1)\},
$$
for all $t$.
\end{lma}

\begin{proof}
Using \fullref{lmanegV}, we deform any given $\pgc$--curve to a generalized $\pgc$--curve with all its vertices at infinity and all curvature concentrations equal to $\pi$ without increasing $\kappa$. In fact, the condition that some curvature concentration is not an integral multiple of $\pi$ implies that $\kappa$ is decreased by this deformation. It is thus sufficient to show that any generalized $\pgc$--curve with all its vertices at infinity and all curvature concentrations equal to $\pi$ may be deformed to a curve in standard position by a generalized $\pgc$--homotopy which increases $\kappa$ arbitrarily little.

Note that any generalized $\pgc$--curve with all its vertices at infinity which has only two infinite vertices is automatically in standard position. Assume inductively that any such curve with $<m$ vertices can be brought to standard position by a $\pgc$--homotopy which increases $\kappa$ arbitrarily little and consider a curve with $m$ vertices.

If the curve has an under-rotated vertex then this vertex can be removed using the method of Case 1 in the proof of \fullref{lmanegV}. The inductive assumption finishes the proof in that case. We thus assume that the curve does not have any under-rotated vertices.

Consider a curve which satisfies this assumption and which has two curvature concentrations of opposite signs. Such a curve must have two vertices $a$ and $b$ with curvature concentrations of opposite signs which are connected by a geodesic arc $C$. Let $A$ and $B$ be the other geodesics with endpoints at $a$ and $b$, respectively. Since every vertex is over-rotated, it follows that $C$ separates $A$ from $B$. Push the second endpoint of $B$ which is not equal to $b$ until it is very close to $a$ and the endpoint of $A$ which is not equal to $a$ until it is very close to $b$. With this done, push the curvature concentrations inwards along $C$ and cancel them increasing $\kappa$ arbitrarily little. The initial phase of the latter homotopy is shown in \fullref{fig:neg5}.

\begin{figure}[ht!]
\begin{center}
\includegraphics[angle=-90, width=8cm]{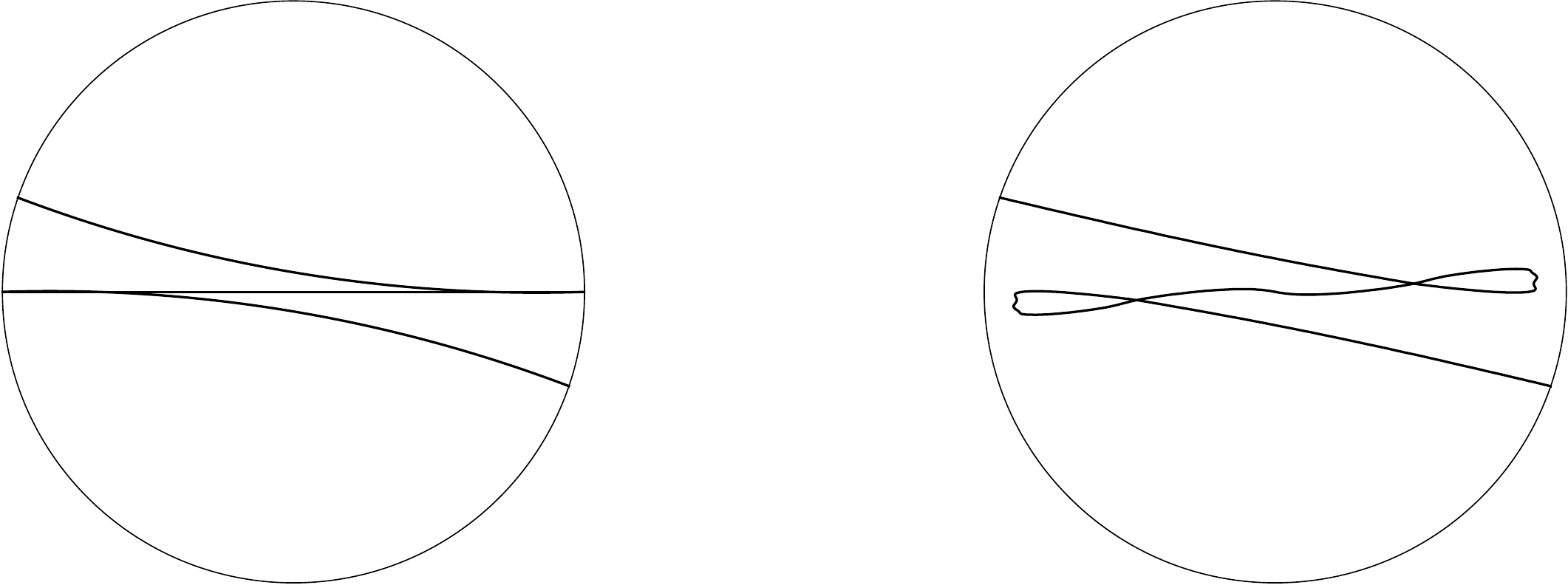}
\end{center}
\caption{Canceling opposite concentrations}
\label{fig:neg5}
\end{figure}

Finally, assume that there are no under-rotated vertices and that all of the vertices have curvature concentrations of the same sign. In this case we pick a first vertex and move the second toward it, eventually creating a small geodesic. We then move the third vertex toward the second and so on. We claim that either this process creates a standard curve or it creates an under-rotated vertex. To see this, let $p_j$ be the first vertex which has not been moved. There are two cases, either $p_j$ can be moved without passing $p_{j+2}$. In this case the construction continues. Or, $p_j$ must pass $p_{j+2}$ in which case an under-rotated vertex at $p_{j+1}$ is created. The claim and the lemma follow.
\end{proof}

\subsection{Constant curvature spheres}\label{ssec:4C}
\begin{lma}\label{lmaposV}
Let $\Sigma$ be the $2$--sphere with a constant curvature metric and let $c$ be a $\pgc$--curve in $\Sigma$. Then there exists a $\pgc$--homotopy $c_t$, $0\le t\le 1$, with $c_0=c$, with $\kappa(c_t)$ a non-increasing function of $t$, and with the following property: the underlying curve of $c_1$ is a geodesic and any curvature concentration of $c_1$ equals $\pi$.
\end{lma}

\begin{proof}
As in \fullref{lmanegV}, scaling invariance of $\kappa$ implies we may assume $K=1$. We claim that any $\pgc$--curve $c$ which does not have underlying curve a geodesic admits a $\pgc$--homotopy which does not increase $\kappa$ and which decreases the number of vertices with curvature concentration not a multiple of $\pi$. As in the proof of \fullref{lmaflatV} this finishes the proof.

Let $p_0$ denote the start and $p_1$ denote the endpoint of a geodesic $e_1$ of $c$. Let $e_0$ and $e_2$ be the other geodesic segments connecting to $p_0$ and $p_1$, respectively, and let $p_2$ denote the other vertex of $e_2$. Let $G_j$ be the great circle in which $e_j$ has its image. If $g$ is a geodesic on $\Sigma$ we let $l(g)$ denote its length. We must consider two separate cases. In the first case $l(e_1)\ne m\pi$ for all integers $m>0$. In this case we deform the curve by moving $p_1$ along $G_2$ in such a way that the length $l(\gamma_0)$ of the vertex curve $\gamma_0$ at $p_0$ {\em decreases}. We stop this deformation the first time $p_1$ hits $p_2$ or one of the points in $G_0\cap G_2$, or, when $l(\gamma_0)=n\pi$, for some integer $n$. At this instant we obtain a curve with the number of vertices with curvature concentration not an integral multiple of $\pi$ one smaller than the corresponding number for $c$.

A straightforward case by case check shows that this deformation does not increase $\kappa$. More precisely, there are 16 subcases to check. They arise as follows. First $l(e_1)=a+2\pi m$, $0<a<\pi$ where $m\ge 0$ is even or odd, second the tangent vectors of $e_0$ and $e_2$ at $p_0$ and $p_1$, respectively, points into different components of $\Sigma-G_1$ or into the same component, third and fourth the lengths of the vertex curve $\gamma_j$ at $p_j$ satisfies $l(\gamma_j)=a+2\pi m$, $0<a<\pi$, where $m\ge 0$ is even or odd, $j=0,1$. However, the fact that $\kappa$ does not increase follows in all of these subcases from one of the following two computations.

First, with notation as in \fullref{fig:pos1}, we calculate
$$
\beta'+(\pi-(\alpha-\alpha'))+\pi-\beta+\int_\Omega K\,dA=2\pi.
$$
Thus,
$$
\Delta\kappa=(\alpha'+\beta')-(\alpha+\beta)=-\int_\Omega K\,dA<0.
$$
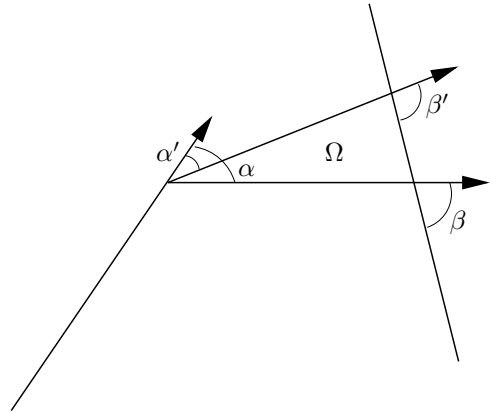
\begin{figure}[ht!]
\begin{center}
\begin{overpic}[angle=-90,width=6.4cm]{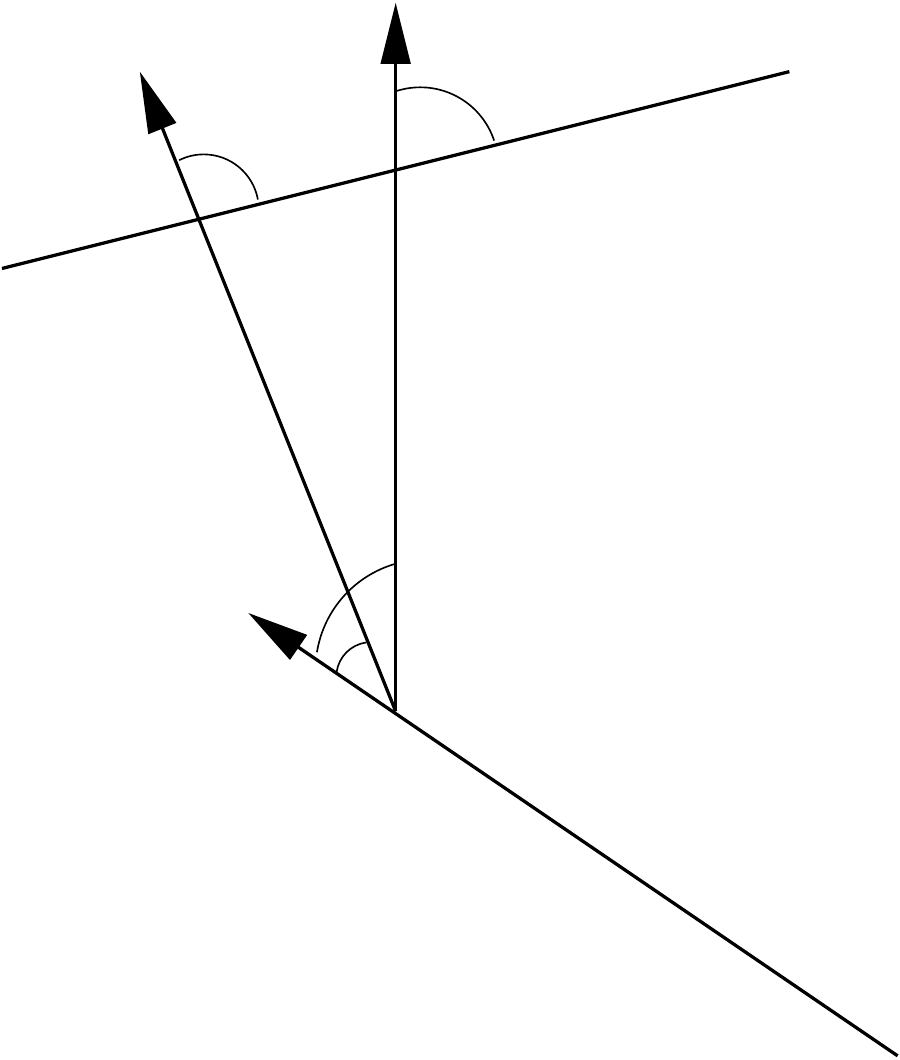}
{\small
\put(47,49){$\alpha$}
\put(30,52){$\alpha'$}
\put(65,52){$\Omega$}
\put(86,62){$\beta'$}
\put(91,38){$\beta$}
}
\end{overpic}
\end{center}
\caption{Removing a vertex I}
\label{fig:pos1}
\end{figure}
Second, with notation as in \fullref{fig:pos2}, we calculate
$$
\beta+(\pi-\beta')+(\pi-(\alpha-\alpha'))+\int_\Omega K\,dA=2\pi.
$$
Thus,
$$
\Delta\kappa=(\alpha'+\beta')-(\alpha+\beta)=
\int_\Omega K\,dA-2(\alpha-\alpha')< 0,
$$
where the last inequality follows since $2(\alpha-\alpha')=\int_\Gamma K\,dA$ where $\Gamma$ is the angular region between the geodesics connecting antipodal points and intersecting at an angle $\alpha-\alpha'$, and since $\Omega\subset\Gamma$.
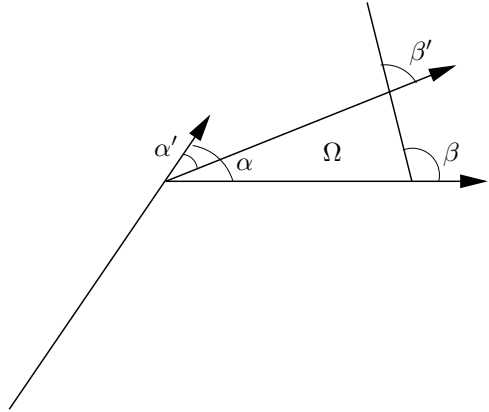
\begin{figure}[ht!]
\begin{center}
\begin{overpic}[angle=-90,width=6.4cm]{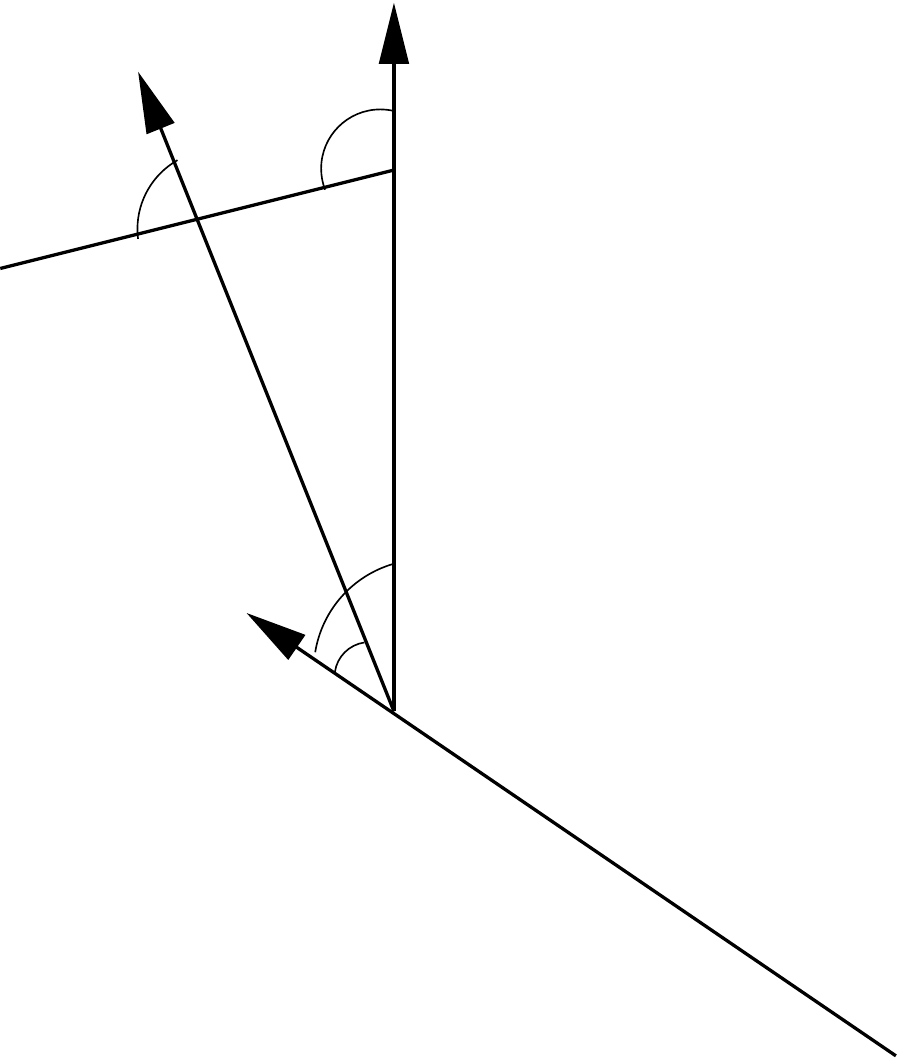}
{\small
\put(47,50){$\alpha$}
\put(30,53){$\alpha'$}
\put(65,52){$\Omega$}
\put(83,73){$\beta'$}
\put(90,52){$\beta$}
}
\end{overpic}\end{center}
\caption{Removing a vertex II}
\label{fig:pos2}
\end{figure}

In the second case $l(e_1)=n \pi$ for some integer $n>0$. In this case we rotate $e_1$ in such a way that the length of at least one of the vertex curves at the endpoints of $e_1$ decreases. The process stops when one of the vertex curves becomes an integral multiple of $\pi$. This finishes the proof.
\end{proof}

\begin{lma}\label{lmapos>}
Let $c_0$ and $c_1$ be any two $\pgc$--homotopic curves on a Riemann surface $\Sigma$ as in \fullref{lmaposV}. Then there exists a $\pgc$--homotopy $c_t$, $0\le t\le 1$, connecting $c_0$ to $c_1$ with
$$
\kappa(c_t)\le\min\{\kappa(c_0),\kappa(c_1),2\pi\},
$$
for all $t$.
\end{lma}

\begin{proof}
\fullref{lmaposV} shows that it is enough to consider two curves with the properties of $c_1$ there. We first show how to deform such a curve to a multiple of a closed geodesic. Fix a first vertex $p$ and the orientation of the great circle of its incoming geodesic. Move the second vertex $q$ to the antipodal point of the first and rotate the arc which connects $p$ and $q$ an angle $\pi$. Note that after the rotation, the orientation of the arc agrees with the fixed one. If the orientation of the vertex curves at $p$ and $q$ are the same then this rotation removes two vertices and decreases $\kappa$. If the two vertex curves have opposite orientations this rotation does not change $\kappa$ and one vertex with vertex curve of length $2\pi$ is created. Splitting this new born vertex into two and repeating the above argument removes them. In this way, we eventually  remove all vertices.

To finish the proof we need only show how to increase the number of times a curve encircles a geodesic by $2$ not increasing $\kappa$ by more than $2\pi$. We use the following procedure. Create two curvature concentrations of length $\pi$ and of opposite orientations keeping $\kappa$ not larger than $2\pi$. Applying the procedure described above, we arrive at a curve which goes two more times around the geodesic.
\end{proof}

\subsection{Curves with self-tangencies on the sphere}\label{ssec:4D}
A generic circle immersion into a surface has only transverse double points. In generic regular homotopies there appear isolated instances of triple points and self-tangencies. The self-tangencies are of two kinds {\em direct} when the tangent vectors at the tangency point agree and {\em opposite} when they do not. Let $\Sigma$ denote the $2$--sphere with a constant curvature metric throughout this subsection.

\begin{lma}\label{lmaost}
Any circle immersion $c\co S^1\to\Sigma$ with an opposite self-tangency satisfies $\kappa(c)>2\pi$.
\end{lma}

\begin{proof}
Since $c$ is not a geodesic we may decrease $\kappa(c)$ by a small deformation, keeping the self-tangency, see the proof of \fullref{prploc->glob}. Let $b$ be the curve resulting from such a deformation. Fix a partition $\pi$ of the circle such that the $\pgc$--approximation $b^\pi$ of $b$ satisfies $\kappa(b^\pi)<\kappa(c)$, see \fullref{lmasmooth->pgc}. Let $\tilde b$ be a $\pgc$--curve close to $b^\pi$ which contains two segments of the same geodesic close to the self-tangency point of $c$ and such that $\kappa(\tilde b)<\kappa(c)$.

Apply the deformation in the proof of \fullref{lmaposV} to $\tilde b$ with the endpoint of one of the self-tangency segments as $p_1$. Note that as this process reaches the endpoint of the other self-tangency segment, the $\pgc$--curve constructed must contain a geodesic segment with at least one vertex curve of length $\pi$. Since the process does not increase $\kappa$ and since the lift of a $\pgc$--curve is closed, it follows that $\kappa(\tilde b)\ge 2\pi$ and therefore $\kappa(c)>2\pi$.
\end{proof}

In order to prove \fullref{thmK=const} (b) we are going to apply
Arnol'd's $J^-$--invariant of immersed curves on $S^2$.  Arnol'd
introduced his invariant for circle immersions in the plane, see
\cite{Ar}. Its existence for curves on more general surfaces was
established by Inshakov \cite{In} and Tchernov \cite{Tch}. The
existence of the $J^-$--invariant stems from the following fact: if
$c_0$ and $c_1$ are two self-transverse regularly homotopic curves on
$S^2$ then the algebraic numbers of opposite self-tangencies in any
two generic regular homotopies $c_t$, $0\le t\le 1$, connecting them
are equal. This number is called the relative $J^-$--invariant of
$c_0$ and $c_1$, we will denote it $\Delta J^-(c_0,c_1)$. To compute
the algebraic number of self-tangencies each self-tangency moment is
equipped with a sign as follows: it is a {\em positive} moment if it
increases the number of double points of the curve, otherwise it is
{\em negative}.

\begin{lma}\label{lmaJ-}
Let $k>0$ be an integer and consider the circle immersions $S^1\to\Sigma$ which go $k$ and $k+2$ times respectively around a geodesic (a great circle). Let $c_0$ and $c_1$ be any small perturbations of these curves. Then $|\Delta J^-(c_0,c_1)|=2$.
\end{lma}

\begin{proof}
Any sufficiently small perturbation of a curve going around a geodesic can be thought of as a multi-graph over the $0$--section in a tubular neighborhood of a great circle. This implies that there exist deformations without opposite self-tangencies connecting any two such perturbations. To compute $\Delta J^-(c_0,c_1)$ it is thus sufficient to pick two perturbations  and count the algebraic number of opposite self-tangencies in one regular homotopy connecting them. In \fullref{fig:J-}, the first picture illustrates a curve in a neighborhood of the north pole which is obtained by shrinking a perturbed multiple of the equator and in the first deformation one kink is pulled over the south pole.
\begin{figure}[ht!]
\begin{center}
\begin{overpic}[angle=-90,width=12cm]{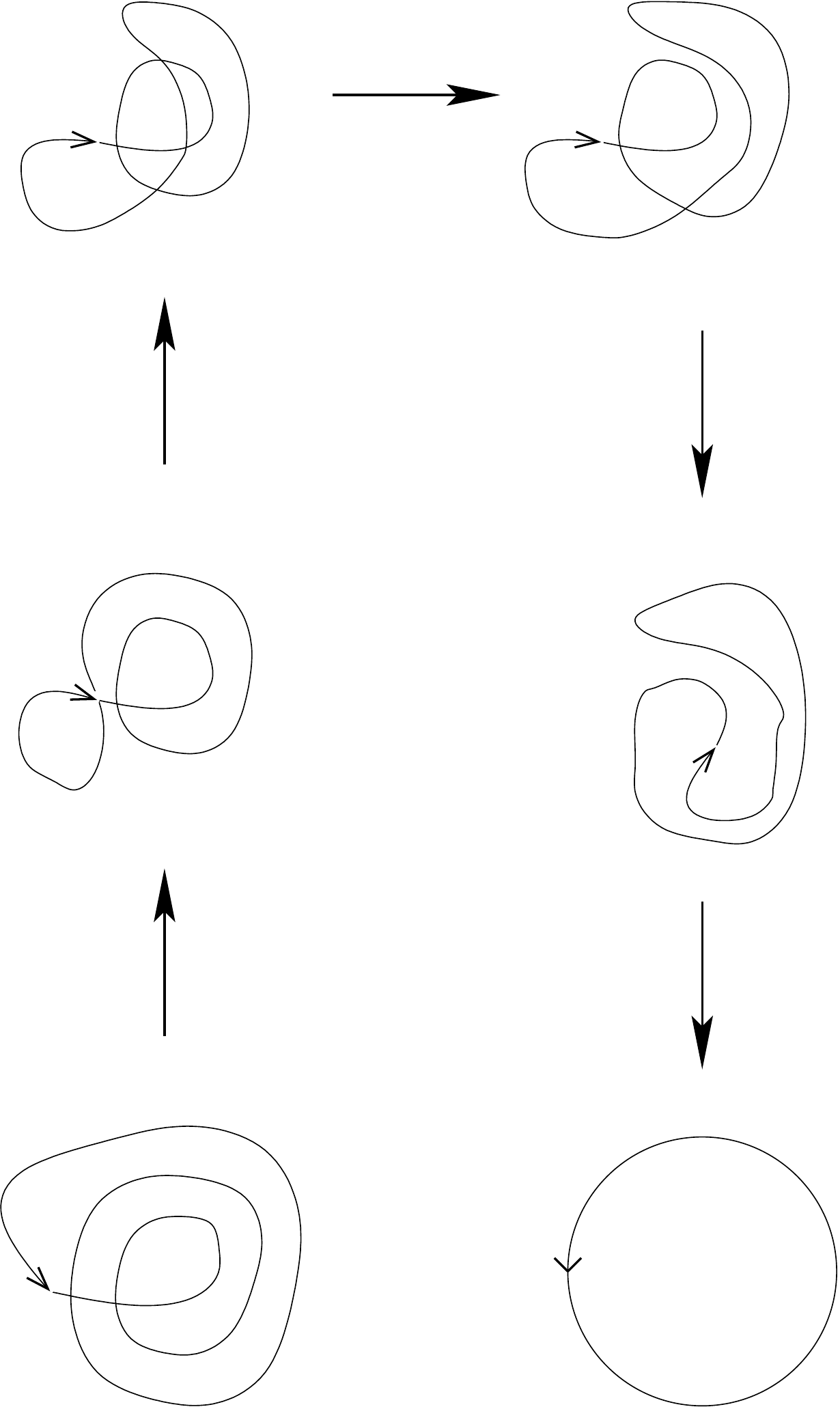}
{\small
\put(70,12){$-1$}
\put(95,30){$-1$}
}
\end{overpic}\end{center}
\caption{Computation of $\Delta J^-$}
\label{fig:J-}
\end{figure}
The lemma follows from \fullref{fig:J-}.
\end{proof}

\subsection{Smoothing $\pgc$--homotopies}\label{ssec:4D1/2}
The $\pgc$--homotopies in \fullref{lmaflatV}, \fullref{lmanegV} and \fullref{lmaposV} have very special forms. We show that such $\pgc$--homotopies can be made smooth in a standard way, increasing the total curvature arbitrarily little.

Consider first two unit vectors $v_{{\rm in}}\in\R^2$ and $v_{{\rm out}}\in\R^2$ and let $\gamma$ be a geodesic in $S^1$ connecting $v_{{\rm in}}$ to $v_{{\rm out}}$ and let $\delta>0$ be given. Fix a family of reference curves $b\bigl(v_{{\rm in}},v_{{\rm out}},\gamma\bigr)$ inside the unit disk such that $b$ agrees with the straight line in direction $v_{{\rm in}}$ ($v_{{\rm out}}$) near the endpoints. The tangent map of $b$ is homotopic to $\gamma$ with endpoints fixed and the total curvature of $b$ exceeds the length of $\gamma$ by at most $\delta$. Clearly there exists such families which depend continuously on the data for all $\delta>0$.

\begin{figure}[ht!]
\begin{center}
\includegraphics[angle=0, width=6cm]{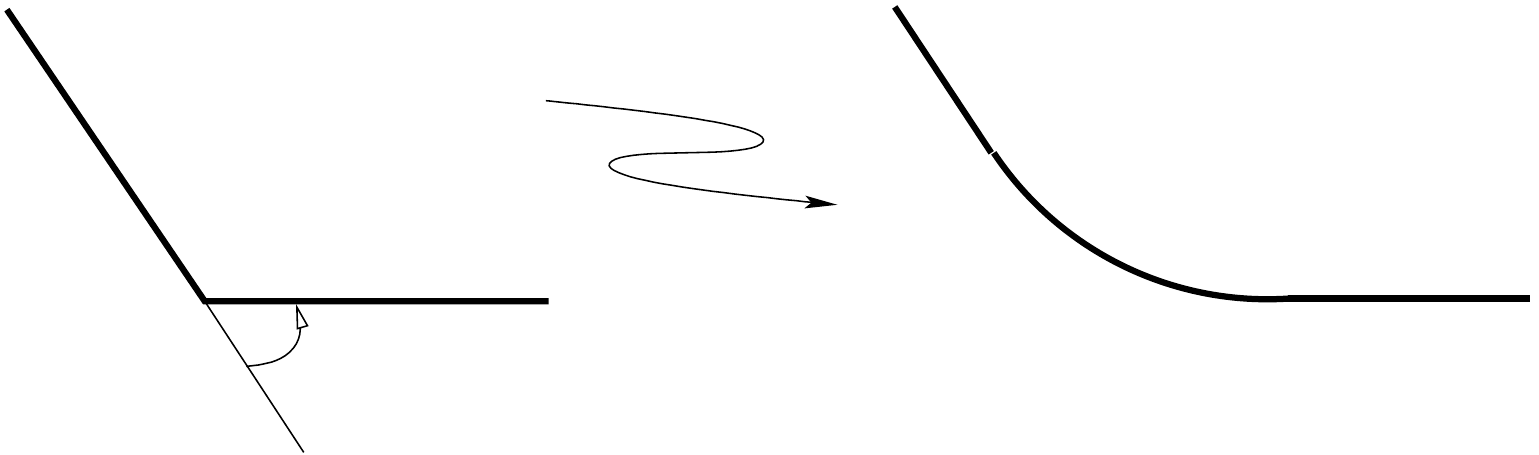}
\end{center}
\caption{Smoothing of exterior angle}
\label{fig:smth1}
\end{figure}

\begin{figure}[ht!]
\begin{center}
\includegraphics[angle=0, width=6cm]{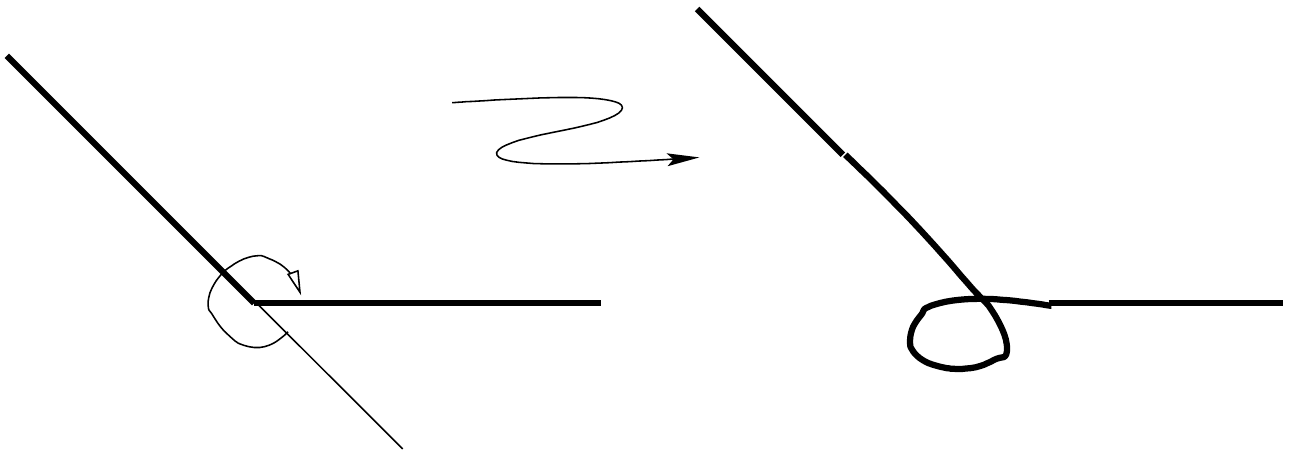}
\end{center}
\caption{Smoothing of interior angle}
\label{fig:smth2}
\end{figure}

We first discuss smoothing of a fixed $\pgc$--curve. Let $c=(c_1,\gamma_1,\dots,c_m,\gamma_m)$ be a $\pgc$--curve such that all vertex curves are immersions. For sufficiently small $\epsilon>0$, we define an $\epsilon$--smoothing of $c$ as follows. Fix a disk of radius $\epsilon>0$ around each vertex of $c$. Under the inverse of the exponential map at a vertex the curve $c$ looks like the model discussed above. More precisely, in the tangent space of the surface $\Sigma$ at a vertex $p$ we have an incoming unit tangent vector $v_{{\rm in}}$ and an outgoing one $v_{{\rm out}}$, and we glue in the curve $b\bigl(v_{{\rm in}},v_{{\rm out}},\gamma\bigr)$ in this tangent space. We then scale the glued in curve by $\epsilon$ and map it back into $\Sigma$ with the exponential map at $p$. Applying this procedure at each vertex we get a smoothing $\tilde c$ of $c$. Note that $\kappa(\tilde c)$ can be made arbitrarily close to $\kappa(c)$ by choosing $\delta$ and $\epsilon$ sufficiently small. (The deviation from the flat case is measured by a curvature integral over a region with area going to $0$ with $\epsilon$.) We illustrate this smoothing process in Figures \ref{fig:smth1} and \ref{fig:smth2}.

We next note that all $\pgc$--homotopies in the proofs of \fullref{lmaflatV}, \fullref{lmanegV} and \fullref{lmaposV} have one of the following forms
\begin{itemize}
\item One geodesic segment of a $\pgc$--curve moves along a consecutive segment.
\item Two curvature concentrations of magnitude $\pi$ and with opposite orientations are created or annihilated somewhere along a $\pgc$--curve.
\item A homotopy of the form presented in the last part of \fullref{lmaflatV}.
\end{itemize}
We call such $\pgc$--homotopies {\em simple}.

\begin{lma}\label{lmasimpler}
Let $c_t$, $0\le t\le 1$, be a simple $\pgc$--homotopy between two $\pgc$--curves $c_0$ and $c_1$ with immersed vertex curves. Then, for any $\epsilon>0$ there exists a regular homotopy $\tilde c_t$ connecting the $\epsilon$--smoothings $\tilde c_0$ and $\tilde c_1$ with
$$
\max_{0\le t\le 1}\kappa(\tilde c_t)\le \max_{0\le t\le 1}\kappa(c_t) + 10\epsilon.
$$
\end{lma}

\begin{proof}
The proof is straightforward. Consider a $\pgc$--homotopy of the first type. We define $\tilde c_t$ as the $\epsilon$--smoothing of $c_t$ for $t$ outside a neighborhood of $1$. Inside a neighborhood of $1$ we may again use local flat models and define $\tilde c_t$ by composing with the exponential map. A picture of such a local model is shown in \fullref{fig:smth3}.
\begin{figure}[ht!]
\begin{center}
\includegraphics[angle=0, width=6cm]{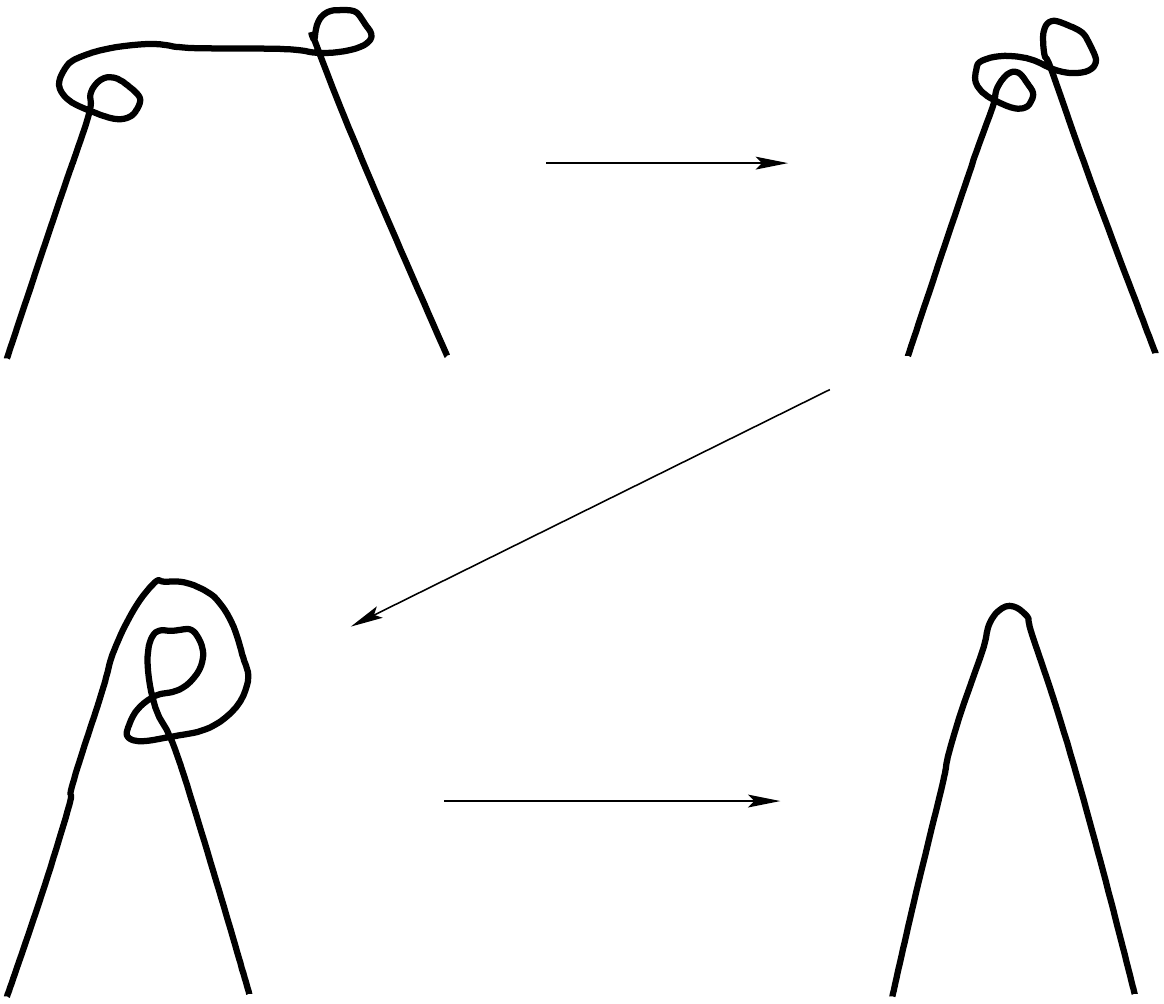}
\end{center}
\caption{A regular homotopy near two meeting vertices}
\label{fig:smth3}
\end{figure}
For the local model of a simple homotopy of the second type, see \fullref{fig:smth4}.
\begin{figure}[ht!]
\begin{center}
\includegraphics[angle=0, width=6cm]{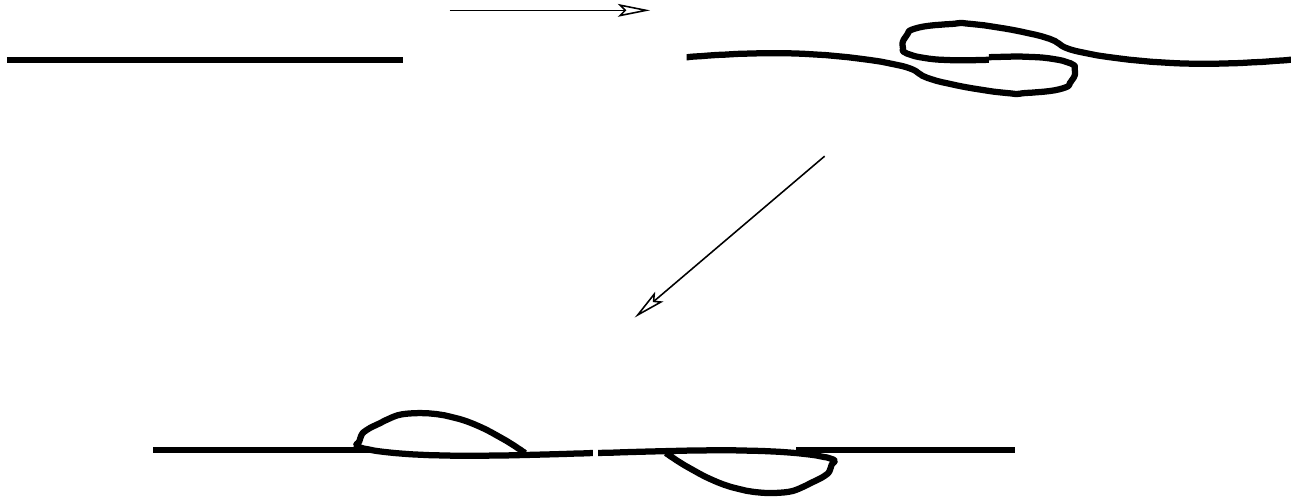}
\end{center}
\caption{Creation of curvature concentrations}
\label{fig:smth4}
\end{figure}

For simple homotopies of the third kind one may simply take $\tilde c_t$ to equal the $\epsilon$--smoothing of $c_t$ for all $t$.
\end{proof}

\subsection{Locally convex curves on flat surfaces}\label{ssec:4E}
Let $\Sigma$ be a flat Riemann surface ($\Sigma\approx T^2$, $\Sigma\approx S^1\times\R$, or $\Sigma\approx\R^2$). For $p\in\Sigma$, $v\in U_p\Sigma$, and $(\xi,m)\in\pi_1(\Sigma)\times\Z$ let $\Omega_{(\xi,m)}(p,v)$ ($\hat\Omega_{(\xi,m)}(p,v)$) denote the space of all (strictly) locally convex circle immersions $c\co S^1\to\Sigma$ of regular homotopy class $(\xi,m)$, see \fullref{sec:1}, such that $c(1)=p$ and $\dot c(1)=v$. (As usual we use the $C^2$--topology on the space of circle immersions.)

\begin{lma}\label{lmaF->*}
The spaces $\Omega_{(\xi,m)}(p,v)$ and $\hat\Omega_{(\xi,m)}(p,v)$ are weakly contractible.
\end{lma}

\begin{proof}
We start in the simply connected case, $\Sigma=\R^2$. Let $F$ denote $\Omega_m(p,v)$ or $\hat\Omega_m(p,v)$. Let $\Gamma\co S^n\to F$ be a continuous map from the $n$--sphere, $\Gamma(x)=c_x\co S^1\to\R^2$. We think of $S^1$ as an interval $[0,L]$ with endpoints identified. Thus $c_x(0)=p$ and $\dot c_x(0)=v$, for some fixed point $p$, some unit vector $v$, and all $x\in S^n$. To prove the lemma we must extend $\Gamma$ continuously to the $(n+1)$--ball $B^{n+1}$.

Fix a small $\epsilon>0$ and a unit vector $w$ such that $\la v, w\ra=0$ and such that in the orientation of the plane induced by the basis $w,v$ the tangential degree of the curves $c_x$ are positive. We claim that there exist continuous maps $t_j\co S^n\to\R$, $j=1,2$, with $0<t_1(x)< t_2(x)< L$, and with the following properties: $\dot c_x(t_1(x))$ lies in the short sub-arc $A_\epsilon$ of $S^1$ between $(\cos\epsilon)v+(\sin\epsilon) w$ and $(\cos 2\epsilon)v+(\sin 2\epsilon) w$, $\dot c_x(t_2(x))$ lies in the short sub-arc $B_\epsilon$ between $(\cos\epsilon)v-(\sin\epsilon) w$ and $(\cos 2\epsilon)v-(\sin 2\epsilon) w$, $t_1(x)$ lies in the component of ${\dot c}_x^{-1}(A_\epsilon)$ closest to $0$, and $t_2(x)$ lies in the component of ${\dot c}_x^{-1}(B_\epsilon)$ closest $L$.

In the strictly locally convex case this claim is obviously true: consider preimages under $\dot c_x$ of fixed points in $A_\epsilon$ and $B_\epsilon$ to define $t_1(x)$ and $t_2(x)$, respectively. To see that it holds also in the non-strictly locally convex case we argue as follows. By continuity, the subset $\Lambda\subset S^n\times S^1$,
$$
\Lambda=\{(x,t)\co \dot c_x(t)\in \inr(A_\epsilon)\},
$$
where $\inr(X)$ denotes the interior of $X$, is open. In particular, for each $x\in S^n$ there exists $r_x>0$ such that
$$
\bigcap_{y\in B(x,r_x)} {\dot c_y}^{-1}(\inr(A_\epsilon))
$$
contains an interval $I_x$. Cover $S^n$ by balls $B(x,r_x)$ with this property. This cover has a Lebesgue number $\delta>0$. Triangulate $S^n$ by simplices which are so small that for every vertex in the triangulation, the union of all simplices in which this vertex lies is a subset of diameter less than $\delta$. It is then straightforward to construct $t_1\co S^n\to\R$: if $v$ is a vertex take $t_1(v)$ as any point in $I_x$ where $x$ is some point such that the union of all simplices containing $v$ lies in $B(x,r_x)$. By contractibility of the interval we can now inductively extend this function over higher dimensional skeleta of the triangulation. In the final step we get the desired function on $S^n$.

Let $q_x$ be the intersection point of the tangent lines to $c_x$ at $c_x(t_1(x))$ and $c_x(t_2(x))$. Define an initial deformation of the curves $c_x$ in the family which pushes $c_x$ toward the piecewise linear curve obtained by replacing $c_x([0,t_1(x)])\cup c_x([t_2(x),L])$ with the curve $\overline{c_x(t_1(x))q_x}\cup \overline{c_x(t_2(x))q_x}$, where $\overline{a b}$ denotes the line segment between $a$ and $b$. To stay in the space of pointed curves we also translate and slightly rotate the curves to ensure that they pass through $p$ with the right tangent, see \fullref{fig:convin}.
\begin{figure}[ht!]
\begin{center}
\begin{overpic}[angle=-90,width=2cm]{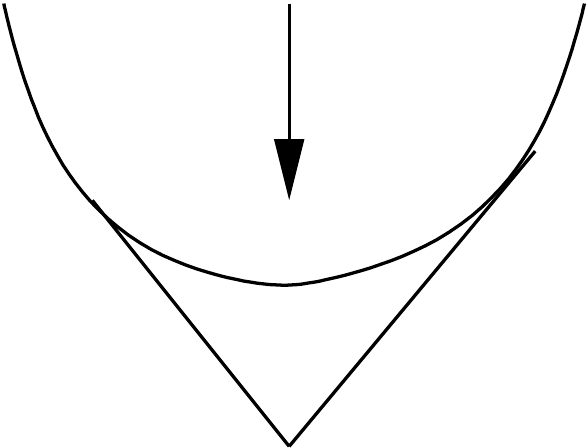}
{\small
\put(80,47){push}
}
\end{overpic}
\end{center}
\caption{An initial deformation}
\label{fig:convin}
\end{figure}
Clearly, this deformation can be made continuous in $x$ and chosen in such a way that the resulting curves, still denoted  $c_x$, are strictly locally convex at $c_x(0)=p$ for each $x\in S^n$.

The next step is to deform the curves so that the marked point is a global minimum of the height function in direction $w$. To this end, let $T_j\co S^n\to\R$, $j=1,2$, $0<T_1(x)<T_2(x)<L$ be continuous functions on $S^n$ such that $\dot c_x(T_1(x))$ ($\dot c_x(T_2(x))$) lies in an $\eta$--arc $C_\eta$ ending at $w$ (beginning at $-w$) in the orientation of $S^1$ determined by $\dot c_x$ and such that $T_1(x)$ ($T_2(x)$) lies in the component of ${\dot c_x}^{-1}(C_\eta)$ closest to $0$ (to $L$). In the strictly locally convex case such functions are easily constructed using suitable points in the preimage ${\dot c_x}^{-1}(\pm w)$. In the non-strictly locally convex case, such functions can be constructed using the same arguments that were used in the construction of the functions $t_j$, $j=1,2$, just given. Let $w_1(x)=\dot c_x(T_1(x))$ and let $w_2(x)=-\dot c_x(T_2(x))$.

Pick $M>0$ such that the minimum point $c_x(q)$ in the $w$--direction of any curve $c_x$, $x\in S^n$, satisfies $\la w, c_x(q)-p\ra>-\frac12 M$. Pick $\eta>0$ sufficiently small so that the intersection point $r$ of the lines in direction $w_1(x)$ and $w_2(x)$ passing through $c_x(T_1(x))$ and $c_x(T_2(x))$ respectively satisfies $\la w, r-c_x(q)\ra<-20 M$. Consider the subdivision of the source circle of $c_x$ into two arcs $D_1(x)$ and $D_2(x)$ as follows: the endpoints of the arcs are $T_1(x)$ and $T_2(x)$, $D_1(x)$ contains the marked point, and $D_2(x)$ does not contain it.

For $0\le s\le \frac12$, define $c_{sx}$ as the curve which consists of the following four pieces: the curve $c_x(D_2(x))$, the line segment $l^1_s$ in direction $w_1(x)$ starting at $c_x(T_1(x))$ and such that the projection of this segment to a line in the $w$--direction has length $10 sM$, a suitably scaled version of the translate of the curve $c_x(D_1(x))$ along $l^1_s$, which is tangent to the line through $c_x(T_2(x))$ in direction $w_2(x)$, and finally the line segment $l^2_s$ in direction $-w_2(x)$ connecting the scaled and translated $c_x(D_1(x))$ to $c_x(D_2(x))$ at $c_x(T_2(x))$, see \fullref{fig:convdi}.
\begin{figure}[ht!]
\begin{center}
\includegraphics[angle=-90, width=8cm]{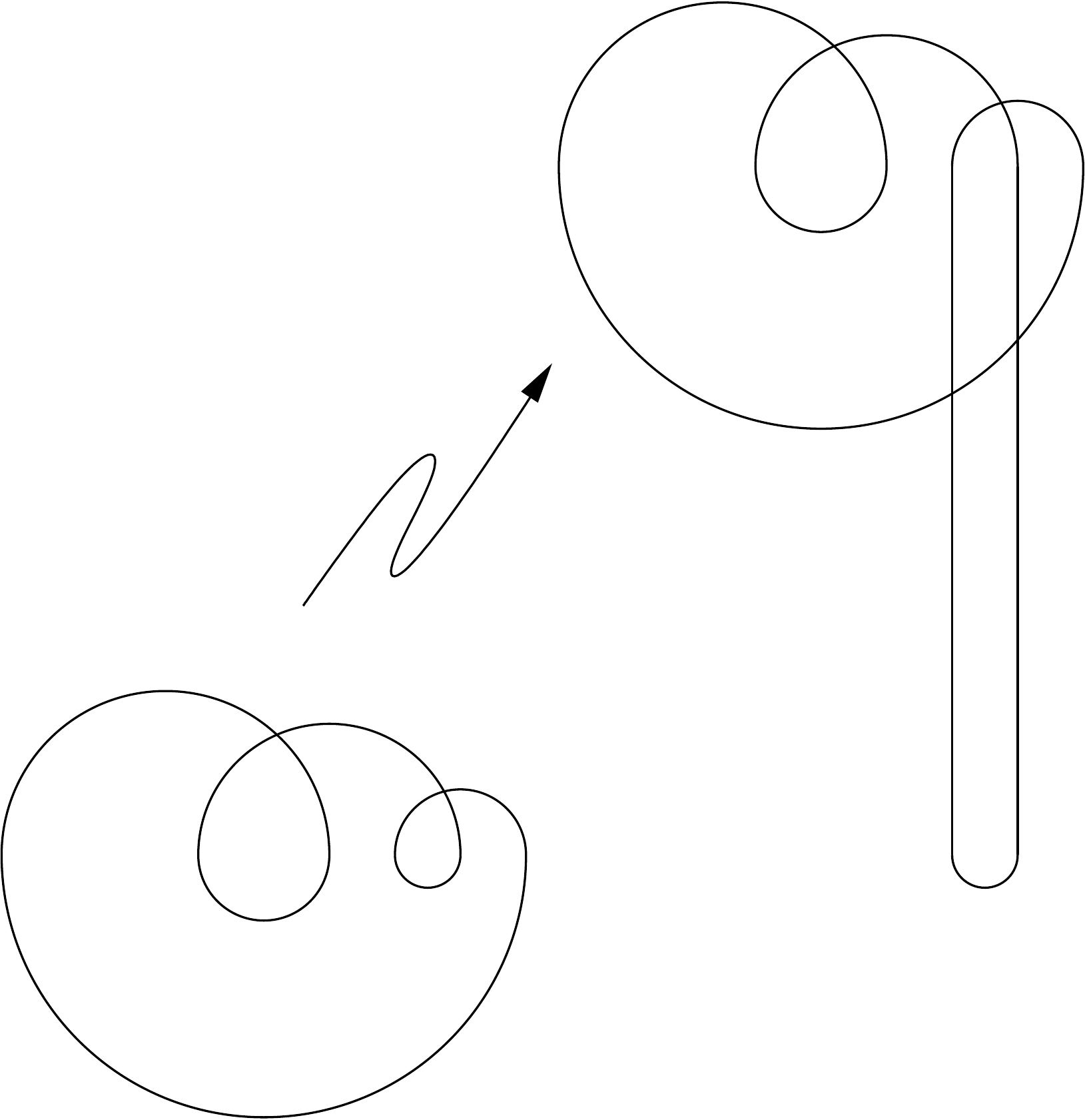}
\end{center}
\caption{Making the marked point a global extremum}
\label{fig:convdi}
\end{figure}
In the case when $F$ is a space of strictly convex curves we replace the straight line segments in $c_{sx}(t)$ above with slightly curved circular arcs. (In fact, since $S^n$ is compact we can take these arcs to have curvature smaller than $\min_{x\in S^n}\min_{t\in[0,L]}(|k_g(t)|)$.) Finally, to have the curves mapping the marked point to $p$ we also compose with a suitable translation.

The second step in the deformation takes all the curves in the family to curves with image in a large circle. Let $C$ be a circle through $p$ with tangent $v$ at $p$. Let $D$ be the bounded component of the plane with boundary $C$. If the radius of $C$ is sufficiently big then we may choose $C$ so that all curves $c_{\frac12 x}$ has image in $D$ and $C\cap c_{\frac12 x}([0,L])=c_{\frac12 x}(0)$. For $0\le s\le 1$ let $c_{(-\frac{s}{2}+\frac12)x}$ be the continuous family of convex curves which is the union of a curve in $C$ starting at $p$ and ending at the intersection point of the negative tangent half-line of $c_{\frac12 x}$ at $c_{\frac12 x}(sL)$, follows this half-line, and then goes along $c_{\frac12 x}$. For $s=1$, we get a curve wrapping around $C$, see \fullref{fig:conv1}.
\begin{figure}[ht!]
\begin{center}
\includegraphics[angle=-90, width=6cm]{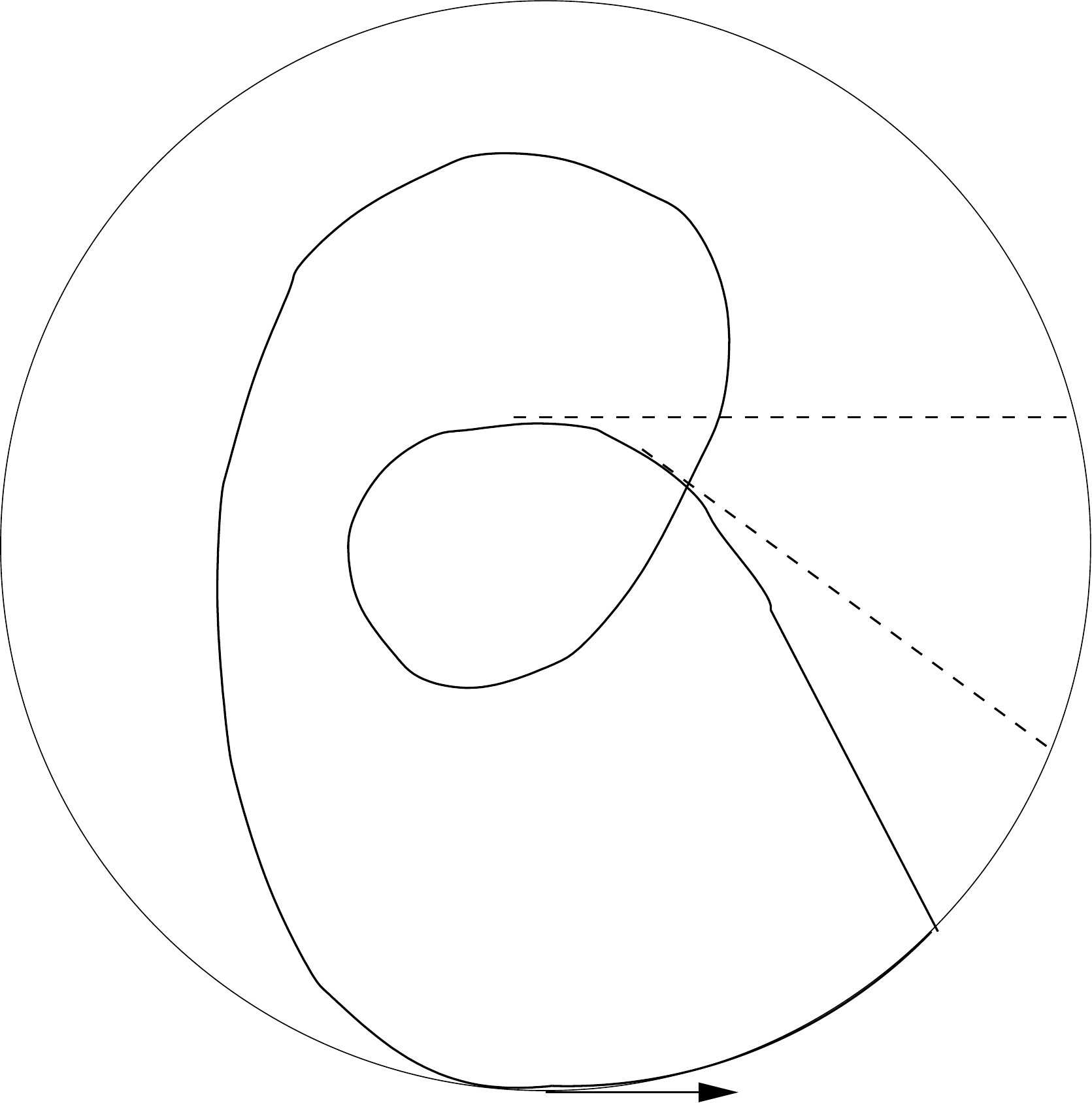}
\end{center}
\caption{Making a locally convex curve circular}
\label{fig:conv1}
\end{figure}
The construction is continuous in $x$ and therefore gives a continuous extension of the map $\Gamma\co S^n\to F$ over $B^{n+1}$, as desired. (Again, in the case of strictly locally convex curves we replace the tangent half-line with a very slightly curved circular arc.)

We next consider the non-simply connected case. In this case we have a covering $\Pi\co \R^2\to\Sigma$ where the deck transformations are translations. We lift the curves $c_x$ at the marked point. In the case when $F$ is the space of non strictly locally convex curves a similar local deformation as the initial deformation above makes all curves strictly locally convex at the marked point. Let $b_x$ denote the lifted curves. We have $b_x(0)=0$, $b_x(L)=q$, $\dot b_x(0)=v=\dot b_x(L)$.

Let $l_0$ and $l_q$ be the straight lines through $0$ and $q$ respectively with tangent vector $v$. By strict local convexity, for one of the unit vectors $w$ orthogonal to $v$, we have $\la b_x(t),w\ra>0$ for $t$ in a (punctured) neighborhood of $0$ and for $t$ in a (punctured) neighborhood of $L$. For convenience, assume that $\la q,w\ra>0$ (otherwise change coordinates in $\R^2$ so that $q=0$). We construct $T_j\co S^n\to\R$, $j=1,2$, analogous to the functions with the same names above, in a similar way as above so that $\dot b_x(T_1(x))\approx w$ and so that $\dot b_x(T_2(x))\approx -w$. As above, we add straight line segments to $b_x$ so that $b_x(0)$ (or $b_x(q)$) is the global minimum of the height function in direction $w$ and so that $b_x$ intersects the region between the two lines $l_0$ and $l_q$ in an arc.

Pick two circles $C_0$ and $C_q$ through $0$ and $q$, respectively, both with tangents $v$ at these points. Let $D_0$ and $D_q$ be the bounded components of the plane which are bounded by $C_0$ and $C_q$, respectively. If the circles have sufficiently large radii then $b_x([0,L])\subset D_0$ for each $x$ and $b_x$ intersects $D_0-D_q$ in an arc. Let $C$ be the boundary of the convex hull of $C_0$ and $C_q$. As above, we deform $b_x$ using its negative tangent half-line. For $0\le s\le\tau$ where $\tau$ is the smallest number such that the negative tangent half-line at $b_x(\tau)$ intersects $C$ in a point in $C_q$. We let $b_{sx}$ be a part of $C$ followed by the negative tangent half-line, in turn followed by the rest of $b_x$. For $s\ge \tau$ we let $b_{sx}$ be the curve with initial part as above followed by a curve in $C_q$, in turn followed by the tangent half-line, and finally the rest of $b_x$. At $s=L$ we find that every curve in the family is a curve which is a segment in $C$ followed by a curve with image in $C_q$, see \fullref{fig:conv2}.
\begin{figure}[ht!]
\begin{center}
\includegraphics[angle=-90, width=4cm]{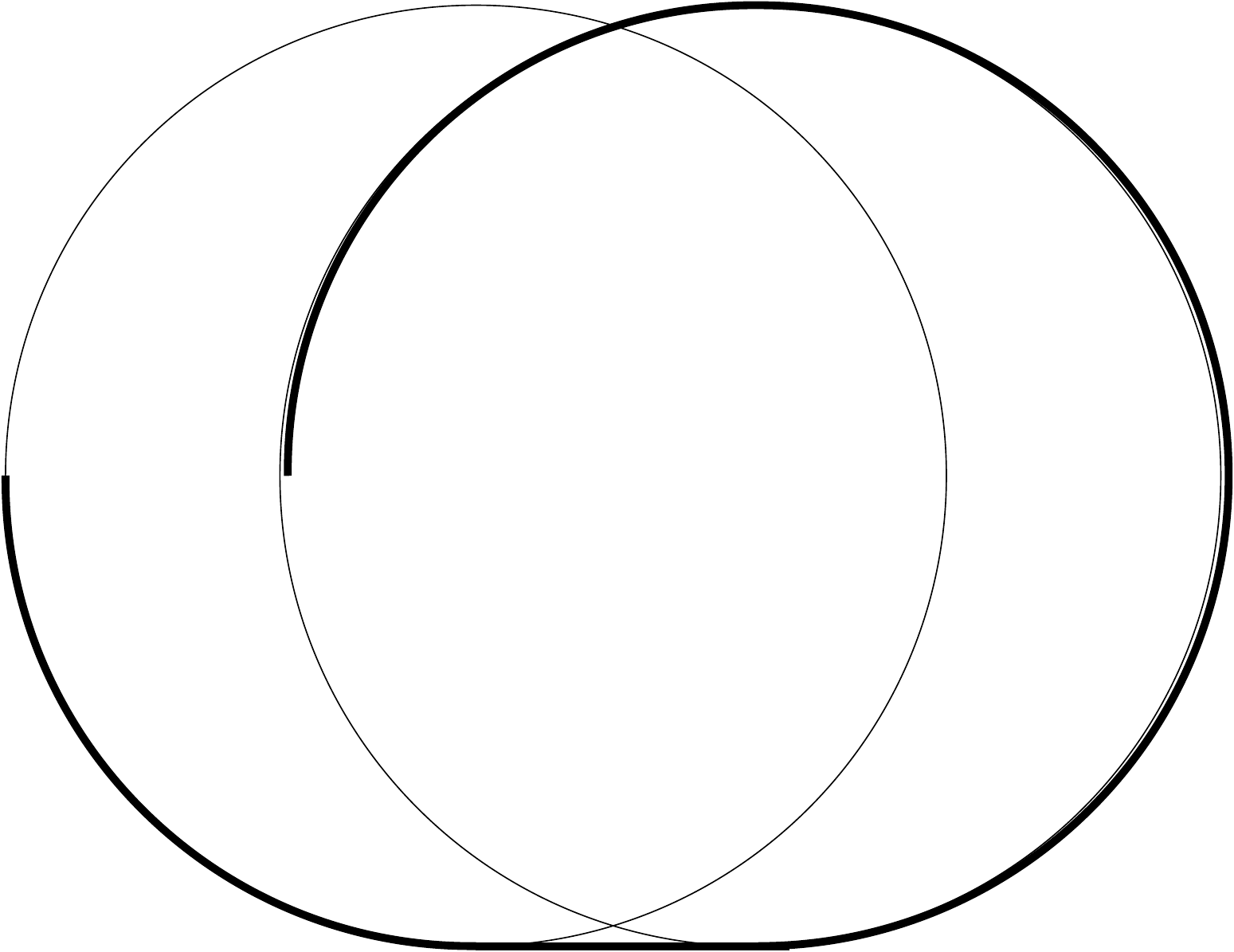}
\end{center}
\caption{The final stage of a deformation of the lift of a locally
convex curve}
\label{fig:conv2}
\end{figure}
In the case of strictly locally convex curves it is easy to modify the deformation described so that it keeps all curves strictly locally convex. This finishes the proof.
\end{proof}

\begin{lma}\label{lmahat=nohat}
The inclusion $\hat\Omega_{(\xi,m)}\subset\Omega_{(\xi,m)}$ induces surjections
$$
\pi_r(\hat\Omega_{(\xi,m)})\to\pi_r(\Omega_{(\xi,m)}),
$$
on homotopy groups, for all $r$.
\end{lma}

\begin{proof}
The proof of \fullref{lmaF->*} homotopes an arbitrary family of (based) locally convex curves to a family of strictly locally convex curves. Hence we need only consider what happens to the base point. Let $\Gamma\co S^n\to\Omega_{(\xi,m)}$ be a family of curves. First use the initial deformation of the proof above to make all curves strictly locally convex in a neighborhood of $1\in S^1$. With this accomplished we lift all curves in the family to $\R^2$, by lifting at $c(1)$. (If $n\ge 2$ then we can lift the whole family in this way however when $n=1$ the start- and endpoints of our lifts may differ by a translation.) Now apply the procedure of the proof of \fullref{lmaF->*}. Note that (in the case $n=1$) the procedure behaves well with respect to translations. Hence, any $\Gamma\co S^n\to\Omega_{(\xi,m)}$ can be homotoped to a map $\Gamma'\co S^n\to\hat\Omega_{(\xi,m)}$. This finishes the proof.
\end{proof}

\section{Proofs}\label{sec:5}
In this section we prove the theorems stated in \fullref{sec:1}.

\subsection{Infima}\label{ssec:5A}
\begin{proof}[Proof of \fullref{thminf}]
Consider first case (c). A well-known theorem of Lyusternik and Fet
\cite{LF} says that the $2$--sphere with any metric (actually any
closed Riemannian manifold) has a non-constant simple closed geodesic,
see also Jost \cite[Section 5.5]{J}. This geodesic traversed once and
twice respectively gives representatives with $\kappa=0$ for both
regular homotopy classes on the $2$--sphere.

In case (b) it is easy to construct a locally convex curve in any regular homotopy class $(\xi,m)$ with $\xi\ne\ast$ or $m\ne 0$. Note that the class $(\ast,0)$ can neither be represented by a closed geodesic nor by a locally convex curve. Hence $\kappa(c)$ can be decreased for each $c$ of regular homotopy class $(\ast,0)$. Moreover, the curvature decreasing procedure in the proof of \fullref{lmaflatV} gives in this case a $\pgc$--curve with underlying curve a geodesic segment. Since the lift of that $\pgc$--curve is closed it must have at least two curvature concentrations of magnitude $\pi$. Hence $\hat\kappa(\ast,0)\ge 2\pi$. Approximating a segment traversed twice with two vertex curves which are rotations $\pi$ and $-\pi$ we find that $\hat\kappa(\ast,0)=2\pi$.

Finally, in case (a),  it follows as above that the infimum is not attained in classes not representable by geodesics. To find the infima in case $\Sigma\approx\R^2$ we apply \fullref{lmanegV} to conclude that it is enough to find the minimal $\kappa$ of a generalized $\pgc$--curve in standard form in a given regular homotopy class. It is straightforward to check that a curve in standard form representing the regular homotopy class $m$ has two vertex curves of length $\pi$ if $|m|=0$ and $|m|+1$ vertex curves of length $\pi$ otherwise.
\end{proof}

\subsection{Locally convex curves}\label{ssec:5B}
\begin{proof}[Proof of \fullref{thmfree}]
Let $v$ be a (covariantly) constant unit vector field on $\Sigma$. Let $\hat\Omega_{(\xi,m)}(v)$ and $\Omega_{(\xi,m)}(v)$ denote the space of strictly locally convex--, respectively, locally convex curves $c$ with $\dot c(1)=v$. Consider the evaluation maps $e(c)=c(1)$,
$$
e\co \hat\Omega_{(\xi,m)}(v)\to\Sigma\quad\text{ and }
e\co \Omega_{(\xi,m)}(v)\to\Sigma.
$$
These maps are clearly Serre fibrations: using translations we can lift any map from an $n$--disk into $\Sigma$. The fibers of these fibrations are $\hat\Omega_{(\xi,m)}(p,v)$ and $\Omega_{(\xi,m)}(p,v)$, respectively, which are both weakly contractible by \fullref{lmaF->*}. Hence
$$
\pi_r(\hat\Omega_{(\xi,m)}(v))\approx\pi_r(\Sigma)
\approx \pi_r(\Omega_{(\xi,m)}(v))
$$
with isomorphisms induced by evaluation.

Assume $m\ne 0$ and consider the fibration
$$
e'\co \hat\Omega_{(\xi,m)}\to S^1,
$$
$e'(c)=\dot c(0)$. This is a Serre fibration since
$\dot c\co S^1\to S^1$ is a covering map
(here it is essential that the curves are {\em strictly} locally
convex). Since the fiber of $e'$ is $\hat\Omega_{(\xi,m)}(v)$ we find that
$$
\pi_r(\hat\Omega_{(\xi,m)})=\pi_r(U\Sigma),
$$
with isomorphism induced by the evaluation map. Since
$\pi_r(U\Sigma)=0$ if $r>1$ we conclude from
\fullref{lmahat=nohat} that
$\pi_k(\Omega_{(\xi,m)})=\pi_k(\hat\Omega_{(\xi,m)})$ for all $k$.

We finally consider $m=0$. In this case the space under consideration is the space of closed geodesics in a fixed homotopy class. If $\Sigma\approx S^1\times\R$ any element in such a space is uniquely determined by its intersection with $\R\times\{1\}$ and $\pi_{S^1}(c(1))$, where $\pi_{S^1}\co S^1\times\R\to S^1$ is the natural projection. Thus $\Omega_{(\xi,0)}\simeq\Sigma$. Similar arguments give the same result when $\Sigma\approx T^2$.
\end{proof}

\subsection{Regular homotopies}\label{ssec:5C}
\begin{proof}[Proof of \fullref{thmK=const}]
Consider first case (a). If $\Sigma$ is flat and both $c_0$ and $c_1$ are closed geodesics or locally convex curves then the theorem follows from \fullref{thmfree}. In all other cases we may first decrease the total curvature a little by \fullref{prploc->glob} and then approximate by a $\pgc$--curve as in \fullref{lmasmooth->pgc}. In the flat case the theorem then follows from \fullref{lmaflatV}, which allows us to deform both curves to a standard form without increasing $\kappa$, \fullref{lmaflat>} which allows us to connect these curves, and \fullref{lmasimpler} which allows us to smooth the entire homotopy keeping control of $\kappa$.  In the case of negative curvature the theorem follows in a similar way from Lemmas \ref{lmanegV}, \ref{lmaneg>} and \ref{lmasimpler}.

In case (b) we argue in the same way to prove the first statement using Lemmas \ref{lmasmooth->pgc}, \ref{lmaposV}, \ref{lmapos>} and  \ref{lmasimpler}. The second statement follows from \fullref{lmaJ-} which shows that any regular homotopy $c_t$, $0\le t\le 1$, connecting the two multiple geodesics of different multiplicities must have an instant $c_{\tau}$ which is a curve with an opposite self-tangency and \fullref{lmaost} which shows that $\kappa(c_{\tau})>2\pi$.
\end{proof}

We end the paper by demonstrating that \fullref{thmK=const} (a) does not hold for Riemann surfaces with metrics of non-constant curvature with $K\le 0$.

\begin{rmk}\label{rmkKle0}
Consider the upper half-plane with coordinates $(x,y)$, $y>0$ and metric
$$
ds^2=e^{2f(y)}(dx^2+dy^2),
$$
where $f\co (0,\infty)\to\R$ is a (weakly) convex function such that its derivative $f'$ satisfies $f'(y)=0$ for $y\in [1,2]\cup[3,4]$, and such that its second derivative $f''$ satisfies $f''(y)>0$ for $y\in (0,1)\cup(4,\infty)$ and such that $f(y)=-\log(y)$ for $y$ in some neighborhood of $0$ and of $\infty$.

Let $A_0=\R\times[1,2]$, $A_1=\R\times[3,4]$, and $B=\R\times(0,\infty)-(A_0\cup A_1)$. Then $K(p)=0$ for $p\in A_0\cup A_1$ and $K(p)<0$ for $p\in B$. Let $c_0$ and $c_1$ be convex curves in $A_0$ and $A_1$, respectively. Then $c_0$ and $c_1$ are regularly homotopic. We claim that for every regular homotopy $c_t$, $0\le t\le 1$ connecting $c_0$ to $c_1$, there exists an instant $c_\tau$ with $\kappa(c_\tau)>2\pi$. To see this note that as long as the curve $c_t$ stays in $A_0$ it must remain convex, otherwise $\kappa> 2\pi$. In particular, the curve must remain embedded. Let $t_0$ be the last moment when the curve lies completely inside the closure of $A_0$. Since embeddedness is an open condition we see that $c_\tau$ is embedded for all $\tau>t_0$ sufficiently close to $t_0$. Note that $c_{t_0}\cap\pa A_0\ne \emptyset$. Pick some line $l$ parallel to $\pa A_0$, and close to a point in $c_{t_0}\cap\pa A_0$ and such that $c_{t_0}$ intersects it transversely in two points. (The existence of such a line follows from Sard's lemma.) Then also $c_\tau$ meets $l$ transversely for $\tau$ sufficiently close to $t_0$. Applying the Gauss--Bonnet theorem to the two curves bounded by the bounded segment of $l$ cut out by $c_\tau$ and the two remaining pieces of $c_\tau$ we find $\kappa(c_\tau)>2\pi$.
\end{rmk}

{\small\parskip 0pt\vskip11pt minus 5pt\relax
{\sl \def\\{\futurelet\next\nocommawithnl}\def\nocommawithnl
  {\ifx\next\newline\else\unskip,\space\ignorespaces\fi}
  \theaddress\par}
{\rightskip0pt plus .4\hsize
{\def\tempab{}\tt\def~{\lower3.5pt\hbox{\char'176}}\def\_{\char'137}%
\ifx\theemail\tempab\else
  \vskip5pt minus 3pt\theemail\par\fi
  \ifx\theurl\tempab\else
  \vskip5pt minus 3pt\theurl\par\fi}
  \vskip11pt minus 5pt
    Received:\qua\receiveddate\ifx\reviseddate\tempab\else \qquad Revised:\qua
  \reviseddate\fi}}
\newpage

\part{Sphere immersions into 3-space}


\count0=493
\setcounter{section}{0}
\title[Regular homotopy and total curvature II]
{Regular homotopy and total curvature II:\\ sphere immersions into $3$--space}

\volumenumber{6}
\issuenumber{}
\publicationyear{2006}
\papernumber{17}
\startpage{493}
\endpage{512}

\doi{}
\MR{}
\Zbl{}

\keyword{immersion}
\keyword{regular homotopy}
\keyword{relatively isotopy tight}
\keyword{sphere eversion}
\keyword{total curvature}

\received{8 February 2005}
\revised{22 February 2006}
\accepted{12 March 2006}
\published{23 March 2006}
\publishedonline{23 March 2006}
\proposed{}
\seconded{}
\corresponding{}
\editor{}
\version{}

\arxivreference{math.GT/0310266}  
\arxivpassword{v7fg8}             



\makeatletter
\def\dnewtheorem#1[#2]#3{\newtheorem{#1}{#3}[section]
\expandafter\let\csname c@#1\endcsname\c@iithm}

\renewcommand{\baselinestretch}{1.04}  

\newtheorem{iithm}{Theorem}[section]
\dnewtheorem{iiprp}[iithm]{Proposition}
\dnewtheorem{iilma}[iithm]{Lemma}

\newtheorem*{iithmKM}{Theorem (Kuiper and Meeks)}

\theoremstyle{definition}
\dnewtheorem{iidfn}[iithm]{Definition}

\makeatother  

\makeautorefname{iiprp}{Proposition}
\makeautorefname{iilma}{Lemma}
\makeautorefname{iithm}{Theorem}

\numberwithin{equation}{section}

\begin{asciiabstract}
We consider properties of the total curvature functional on the space
of 2-sphere immersions into 3-space. We show that the infimum over all
sphere eversions of the maximum of the total curvature during an
eversion is at most 8pi and we establish a non-injectivity result for
local minima.
\end{asciiabstract}

\begin{htmlabstract}
We consider properties of the total curvature functional on the space
of 2&ndash;sphere immersions into 3&ndash;space. We show that the infimum
over all sphere eversions of the maximum of the total curvature during
an eversion is at most 8&pi; and we establish a non-injectivity
result for local minima.
\end{htmlabstract}

\begin{abstract}
We consider properties of the total curvature functional on the space
of 2--sphere immersions into 3--space. We show that the infimum
over all sphere eversions of the maximum of the total curvature during
an eversion is at most $8\pi$ and we establish a non-injectivity
result for local minima.
\end{abstract}

\maketitle

\section{Introduction}\label{iisec:1}
An {\em immersion} of manifolds is a map with everywhere injective
differential. Two immersions are {\em regularly homotopic} if there
exists a continuous $1$--parameter family of immersions connecting one
to the other. The Smale--Hirsch $h$--principle \cite{iiSm2,iiHi} says that
the space of immersions $M\to N$, $\dim(M)<\dim(N)$ is homotopy
equivalent to the space of injective bundle maps $TM\to TN$. In
contrast to differential topological properties, differential
geometric properties of immersions do not in general satisfy
$h$--principles, see \cite[(A) on page 62]{iiGr}.  In this paper and
the predecessor \cite{iiEI}, we study some aspects of the differential
geometry of immersions and regular homotopies in the most basic cases
of codimension one immersions. We investigate whether or not it is
possible to perform topological constructions while keeping control of
certain geometric quantities.

Consider immersions $S^2\to\R^3$. In this case the $h$--principle
implies a famous theorem of Smale \cite{iiSm1}: all immersions
$S^2\to\R^3$ are regularly homotopic. In particular, there exists {\em
sphere eversions} (ie, regular homotopies connecting the unit
$2$--sphere in $\R^3$ to the same immersion with the opposite
(co)orientation).

The {\em total curvature} $\kappa(f)$ of a sphere immersion $f\co S^2\to\R^3$ is the mapping area of its Gauss map:
\begin{equation*}
\kappa(f)=\int_{S^2} |K|\,d\sigma,
\end{equation*}
where $K$ is the Gaussian curvature and $d\sigma$ the area element
induced on $S^2$ by the immersion $f$. The functional $\kappa$ (and
its higher dimensional generalizations) is also known as the
Lipschitz--Killing curvature. It follows from the Gauss--Bonnet
theorem that the global minimum of $\kappa$ equals $4\pi$. An
immersion for which this value is attained is called tight (or convex)
and is known to be the boundary of a convex body in $\R^3$, Kuiper
\cite{iiKu}. This characterization of tight immersions immediately
implies that any sphere eversion $f_t$, $0\le t\le 1$, has an instant
$f_{\tau}$ with $\kappa(f_{\tau})>4\pi$. The following result gives a
corresponding upper bound.

\begin{iithm}\label{thmS^2ev}
For every $\epsilon>0$ there exists a sphere eversion $f_t$, $0\le t\le 1$, with
\begin{equation*}
\max_{0\le t\le 1}\kappa(f_t)<8\pi+\epsilon.
\end{equation*}
In fact, $f_t$ can be chosen so that for each $t\in[0,1]$ there exists a unit vector $v_t\in\R^3$ such that the height function $\la f_t,v_t\ra\co S^2\to\R$ has exactly two non-degenerate critical points.
\end{iithm}
\fullref{thmS^2ev} is proved in \fullref{ssec:6B}. It is
interesting to compare this result to the fact that the $L^2$--norm of
the mean curvature $W(f)$ (also known as the Willmore energy) of
immersions $f\co S^2\to\R^3$ has the following property. For any
sphere eversion $f_t$, $0\le t\le 1$, there exists some $\tau\in(0,1)$
such that $W(f_\tau)\ge 16\pi$. (This is a consequence of two results:
any sphere eversion has a quadruple point, Max--Banchoff \cite{iiBaMa}, and
$W(f)\ge 16\pi$ for any immersion $f$ with a quadruple point, Li--Yau
\cite{iiLiYa}.)

It is unknown to the author whether the result in \fullref{thmS^2ev} is best possible. We therefore  ask: {\em what is the infimum of $\max_{0\le t\le 1}\kappa(f_t)$ over all sphere eversions $f_t$, $0\le t\le 1$?}

In Propositions \ref{prpGauss1} and \ref{prpGauss2} we make two general
observations about local minima of $\kappa$: an immersion of a closed
$n$--manifold into $\R^{n+1}$ which is a critical point of $\kappa$ must
have total curvature a multiple of the volume of the unit $n$--sphere and
its Gauss map cannot have fold singularities. These observations imply
that all local minima of $\kappa$ with curvature function which meets a
non-degeneracy condition are so called {\em relatively isotopy tight
(RIT)} immersions with certain special properties, see \fullref{ssec:7E}.
(Relatively isotopy tight immersions were introduced by Kuiper and Meeks
in \cite{iiKuMe}, we recall their definition in \fullref{ssec:7D}.) We show in \fullref{thmRIT} that RIT immersions with the special properties just mentioned are non-injective. As a consequence we obtain the following result.

\begin{iithm}\label{thmgenlocmin}
Any local minimum $f\co S^2\to\R^3$ of $\kappa$ with $0$ a non-trivial regular value of $K$ is an RIT immersion such that any component of $K^{-1}((-\infty,0])$ is an annulus. In particular, no such local minimum is an embedding.
\end{iithm}
\fullref{thmgenlocmin} is proved in \fullref{ssec:7F}.

\subsection*{Acknowledgements}
During the preparation of this paper, the author was a research fellow
of the Royal Swedish Academy of Sciences sponsored by the Knut and
Alice Wallenberg foundation. He also acknowledges support from the
Alfred P Sloan Foundation and from NSF--grant DMS-0505076.

\section{Total curvature and a sphere eversion}\label{sec:6}
In this section we recall properties of the total curvature of a closed $n$--manifold immersed into $\R^{n+1}$. We then describe a sphere eversion which constitutes a proof of \fullref{thmS^2ev}.

\subsection{Total curvature}\label{ssec:6A}
Let $f\co M\to\R^{n+1}$ be an immersion of a closed oriented
$n$--dimensional manifold. Let $\nu\co M\to\S^n$, where $\S^n$ is the unit
sphere in $\R^{n+1}$, be its Gauss map which takes a point $p\in M$ to the
positive unit normal of $df(T_pM)$. The {\em total curvature} $\kappa(f)$
of $f$ is the mapping area of $\nu$. Standard results, see Cecil and Ryan
\cite{iiCeRy}, then imply
$$
\kappa(f)=\int_M|J(\nu)|\,d\sigma=\frac12\int_{\S^n}\mu(f_v)\,dv,
$$
where in the first expression $|J(\nu)|$ is the absolute value of the Gauss--Kronecker curvature of $f$ and $d\sigma$ is the volume element of the metric induced by $f$, and where in the second, integral geometric expression, $\mu(f_v)$ denotes the number of critical points of the function $f_v\co M\to\R$, $f_v(x)=\la f(x),v\ra$, where $v\in\S^n$ is a unit vector.

\subsection{A sphere eversion}\label{ssec:6B}
Let $f\co S^2\to\R^3$ be any immersion. Then for generic unit vectors $v$ in $\R^3$ the composition of $f$ with the orthogonal projection $\pi_v$ along $v$ to a plane perpendicular to $v$ is a map with only stable singularities. In particular this means that the image of the singular set of $\pi_v\circ f$ is a piecewise regular curve with cusp singularities. More precisely, the jet transversality theorem implies that for generic $v\in\S^2$ the set of points $p\in S^2$ such that $v\in df(T_p S^2)$ is a smooth $1$--dimensional submanifold $\Sigma_v(f)\subset S^2$ and that $v$ gives a vector field along $\Sigma_v(f)$ which is tangent to $\Sigma_v(f)$ with order one tangencies at isolated points in $\Sigma_v(f)$. These points form a subset which we denote $\Sigma'_v(f)\subset \Sigma_v(f)$. We call the image of the singular set $\Sigma_v(f)$ under $\pi_v\circ f$ {\em the fold curve} of $\pi_v\circ f$. It is a piecewise regular curve with singularities corresponding to the points in $\Sigma'_v(f)$, where the fold curve looks like a semi-cubical cusp $\{(x,y)\in\R^2\co x^2=y^3\}$ up to left--right action of diffeomorphisms of the plane, we call such singularities {\em cusps}. Define the {\em total curvature} $\kappa$ of a piecewise regular planar curve with cusps as the sum of the total curvatures $\int_c|k|\,ds$ of its regular pieces $c$, where $k$ is the curvature function of the regular curve $c$. (Note that there are no curvature concentrations at the cusps.) We then have the following
\begin{iilma}\label{lmafoldcurv}
If $f\co S^2\to\R^3$ is an immersion then
\begin{equation}\label{eqintgeom}
\kappa(f)=\frac{1}{2\pi}\int_{v\in \S^2}\kappa(\pi_v\circ f(\Sigma_v(f)))\,dv.
\end{equation}
\end{iilma}
\begin{proof}
This is immediate from the integral geometric expression for the total curvature: the local extrema of  the height function $f_w$ for $w$ perpendicular to $v$ are in one to one correspondence with the local extrema of the height function in direction $w$ on the curve $\pi\circ f(\Sigma_v(f))$ which are not cusps, see \eqref{eqnorform} below. This implies that \eqref{eqintgeom} holds up to an over all normalizing constant. Considering the round sphere it is easy to see that this constant equals $\frac{1}{2\pi}$ as claimed.
\end{proof}

Consider an immersion $f\co S^2\to\R^3$ and let $v$ be a generic unit vector. Thinking of $\R^3$ as $\R^2\times\R$, where $v$ points in the $\R$--direction we write $f=(\pi_v\circ f,f_v)$. For $\lambda>0$ the map $f^{\lambda,v}=(\pi_v\circ f,\lambda f_v)$ is an immersion.
\begin{iilma}\label{lmacountmax}
As $\lambda\to 0$, $\kappa(f^{\lambda,v})$ tends to twice the total curvature of the fold curve of $\pi_v\circ f$. As $\lambda\to\infty$, $\kappa(f^{\lambda,v})$ tends to $2\pi$ times the number of critical points of the Morse function $f_v$.
\end{iilma}
\begin{proof}
The first statement is a consequence of \fullref{lmafoldcurv} together with the fact if $w\ne \pm v$ is any unit vector, $w=\alpha v+\beta w_0$, $\la w_0,v\ra=0$, $\la w_0,w_0\ra=1$, and $\beta>0$, then ${f^{\lambda,v}}_w\to f_{w_0}$ as $\lambda\to 0$. The second statement follows similarly from the integral geometric formula for the total curvature.
\end{proof}

\begin{figure}[htbp]
\begin{center}
\begin{overpic}[angle=-90,width=7cm]{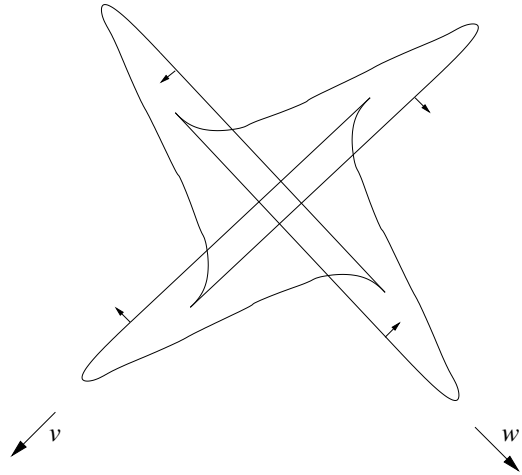}
{\small
\put(7,6){$v$}
\put(93,6){$w$}
}
\end{overpic}
\end{center}
\caption{Middle stage of a sphere eversion}
\label{fig:inout1}
\end{figure}

\begin{figure}[htbp]
\begin{center}
\includegraphics[angle=-90, width=10cm]{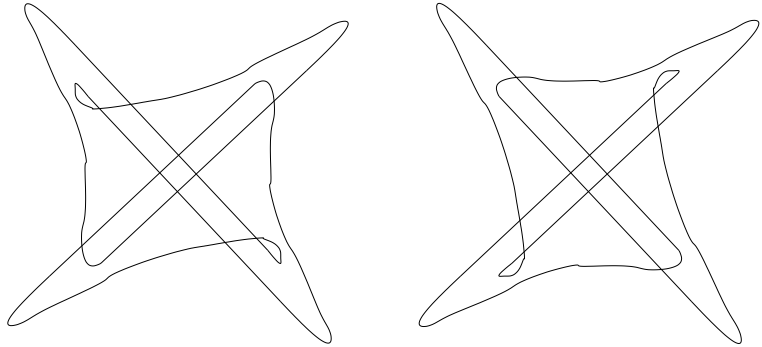}
\end{center}
\caption{Two immersed disks}
\label{fig:inout2}
\end{figure}

We will next construct a sphere eversion $f_t\co S^2\to\R^3$, $0\le t\le 1$, such that\break $\max_{0\le t\le 1}\kappa(f_t)\le 8\pi +\epsilon$. To this end we first describe the middle stage $g_{\frac12}$ of a sphere eversion $g_t\co S^2\to\R^3$, $0\le t\le 1$, closely related to $f_t$. This middle stage is a nearly planar immersion: $g_{\frac12}=(\pi_v\circ h,\lambda h_v)$ for some immersion $h$ and for small $\lambda>0$. In \fullref{fig:inout1} we depict the image of the projection $\pi_v\circ h$. 
In order to see that \fullref{fig:inout1} is really the projection of a sphere immersion we argue as follows. Subdivide $S^2$ into three parts: an thin band $B$ around the equator, a northern disk $D_n$, and a southern disk $D_s$. \fullref{fig:inout2} depicts the images under $\pi_v\circ h$ of $D_n$ and $D_s$. In order to connect these pieces with the band $B$ to form an immersed sphere we consider the local model of an immersion projecting with a cusp. If $(u,v)$ are local coordinates around a point in $S^2$ and $(x,y,z)$ are coordinates in $\R^3$ we have
\begin{align}\notag
x(u,v)&=u,\\\notag
y(u,v)&=v(v^2+u),\\\label{eqnorform}
z(u,v)&=v.
\end{align}
The fold curve of the projection to the $xy$--plane is given by $\{(u,v)\co u=-3v^2\}$. The projection of the curve $\{(u,v)\co u=-3v^2+\epsilon\}$ to the $xy$--plane is a curve with a kink and one transverse double point for $\epsilon>0$, as $\epsilon\to 0$ the kink shrinks, for $\epsilon=0$ the curve has a cusp, and for $\epsilon<0$ the curve is injective. Using this model around the small kinks of the curves $\pi_v\circ h(\pa D_n)$ and $\pi_v\circ h (\pa D_s)$, see \fullref{fig:inout2}, it is straightforward to connect $h(D_n)$ and $h(D_s)$ with an immersion $h$ of $B$ so that $\pi_v\circ h$ is as in \fullref{fig:inout1}.

In order to describe the sphere eversion $g_t$ it thus remains to connect the middle stage $g_{\frac12}$ described above with standard spheres of opposite coorientations. To this end, we use an abstract argument (although it is straightforward to draw a sequence of somewhat complicated pictures).

Consider any sphere immersion $f\co S^2\to\R^3$ such that the height function $f_v(p)=\la v,f(p)\ra$ has exactly two non-degenerate critical points for some direction $v$. We say that such a critical point is {\em positive (negative)} if the direction of the coorienting normal of the immersion agrees (does not agree) with its mean curvature vector at that point. Let the minimum of $f_v$ be  $f_v(q_-)=m$ and the maximum be $f_v(q_+)=M$. Then for $m<a<M$, $f(f_v^{-1}(a))$ is an immersed planar curve in a plane perpendicular to $v$. In this way we can view $f$ as a $1$--parameter family of immersed planar curves which begins and ends at embedded circles. In particular, since the tangential degree of the members of such a family is constant, it follows that the two critical points in the direction $v$ have the same sign.

\begin{iilma}\label{lmatwoextr}
Let $f_j\co S^2\to\R^3$ be immersions and let $v_j$ be unit vectors, $j=0,1$. Assume that the height functions $(f_j)_{v_j}$, $j=0,1$, have exactly two non-degenerate critical points and that the signs of the critical points of $(f_0)_{v_0}$ agree with the signs of the critical points of $(f_1)_{v_1}$. Then there exists a regular homotopy $f_t\co S^2\to\R^3$, $0\le t\le 1$, from $f_0$ to $f_1$, and a continuous family $v_t$, $0\le t\le 1$, of unit vectors such that $(f_t)_{v_t}$ has exactly two non-degenerate critical points for all $t$.
\end{iilma}

\begin{proof}
After composing $f_0$ with a rotation of $\R^3$ which takes $v_0$ to
$v_1$, and a scaling and translation in the $v_1$--direction we may
assume that $v_0=v_1=v$ and that the maxima (minima) of the functions
$(f_0)_v$ and $(f_1)_v$ both equal $M$ ($m$). The lemma is then a
consequence of the Smale--Hirsch $h$--principle for immersed planar
circles as follows. We view the two immersions as two paths $F_0\co
S^1\times[0,1]\to\R^2$ and $F_1\co S^1\times[0,1]\to\R^2$ of immersed
plane curves, where the curve $F_j(\bullet,\tau)$ is
$f_j\bigl((f_j)_v^{-1}(m+\delta+\tau(M-m-2\delta))\bigr)$, $j=0,1$,
for some small $\delta>0$. Then the start- and end-curves
$F_j(\bullet,0)$ and $F_j(\bullet,1)$, $j=0,1$, are very close to
simple convex planar curves since they are intersections of the
original immersions with planes perpendicular to the direction of the
height function very close to its extrema. The Smale--Hirsch
$h$--principle says that these paths of curves can be extended to a
family of plane curve immersions $F\co S^1\times[0,1]^2\to\R^2$ with
$F(\bullet,\bullet,j)=F_j$, $j=0,1$, if and only if the corresponding
paths $\vec F_j\co S^1\times[0,1]\to\R^2-\{0\}$, which assign to a
point $(\theta,t)\in S^1\times[0,1]$ the non-zero tangent vector
$\frac{\pa}{\pa\theta} F_j(\theta,t)$, can be similarly extended as a
continuous map $\vec F$ into $\R^2-\{0\}$.

Such a family $F\co S^1\times[0,1]^2\to\R^2$ of plane curves can be used in an obvious manner to construct a regular homotopy $f_t\co S^2\to\R^3$, $0\le t\le 1$, which connects $f_0$ to $f_1$ and which has properties as claimed. In order to finish the proof of the lemma we thus need only show that the topological extension problem has a solution. This is straightforward: consider $S^1$ as a union of a $0$--cell $e_0$ and a $1$--cell $e_1\approx[0,1]$ with both endpoints identified with $e_0$. First define an extension of $\vec F$ on $e_0\times\{j\}\times[0,1]$ connecting $\vec F_0(e_0,j)$ to $\vec F_1(e_0,j)$, $j=0,1$, in such a way that the loop $\vec F$ on $e_0\times\pa[0,1]^2$ is contractible. Then extend it to $\vec F\co e_0\times[0,1]^2\to\R^2-\{0\}$. Since, $\vec F_0|e_1\times\{j\}$ is homotopic to $\vec F_1|e_1\times\{j\}$, $j=0,1$, we can define $\vec F$ also on $e_1\times\{j\}\times[0,1]$, $j=0,1$. With this done $\vec F$ is defined on the boundary of $[0,1]^3$. Since $\pi_2(\R^2-\{0\})=0$ we can find the desired $\vec F$ by extension over the cube.
\end{proof}

With \fullref{lmatwoextr} established, we return to our sphere eversion $g_t$. If the immersion $g_{\frac12}$ depicted in \fullref{fig:inout1} is sufficiently elongated in the $v$--direction then the resulting immersion has exactly two critical points in the $v$--direction which are both negative. \fullref{lmatwoextr} then implies that we can connect $g_{\frac12}$ with the standard immersion with the outward coorientation. On the other hand, if $g_{\frac12}$ is sufficiently elongated in the $w$--direction then the resulting immersion has exactly two critical points in the $w$--direction which are both positive and \fullref{lmatwoextr} implies that $g_{\frac12}$ can be connected with the standard immersion with the inward coorientation as well. This gives a sphere eversion $g_t$, $0\le t\le 1$.

\begin{proof}[Proof of \fullref{thmS^2ev}]
We show that the sphere eversion $g_t$ described above can be carried out in such a way that $\kappa$ exceeds $8\pi$ by an arbitrarily small amount. As shown in \fullref{fig:inout25}, we can deform the middle stage depicted in \fullref{fig:inout1} keeping it an immersion until it is arbitrarily close to the degenerate middle stage depicted on the right in \fullref{fig:inout25}. Furthermore, as the middle stage approaches the degenerate middle stage, it is clear from the picture that the total curvature of the fold curve approaches the sum of the angles marked on the degenerate middle stage (exterior angles at ordinary corners of the fold curve, interior angles at cusps).
\begin{figure}[htbp]
\begin{center}
\includegraphics[angle=0, width=9cm]{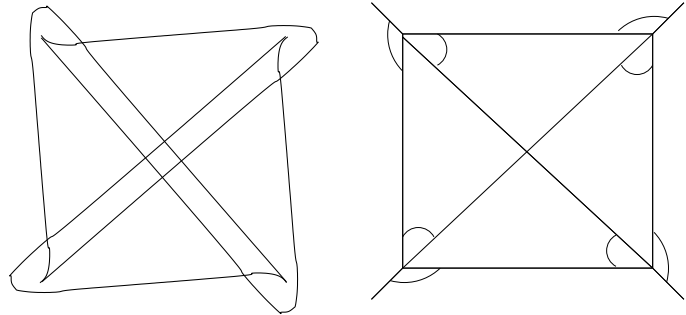}
\end{center}
\caption{Deforming towards a degenerate middle stage}
\label{fig:inout25}
\end{figure}
In \fullref{fig:inout3}, we show how to connect the degenerate middle stage to an immersion which is very elongated in one direction in which it has only two local extrema. The angles contributing to the total curvature of the fold line of the degenerate immersion are easily seen to have sum equal to $4\pi$ at the middle stage and to decrease toward $2\pi$ as we rotate to obtain the elongated version on the left in \fullref{fig:inout3}. It follows that we can connect two immersions $g_0$ and $g_1$ with exactly two extrema in directions $v$ and $w$, respectively, where $\la v,w\ra=0$, which are very elongated in their respective directions and such that the critical points of $g_0$ in direction $v$ have signs opposite to those of $g_1$ in direction $w$, by a regular homotopy $g_t$, $0\le t\le 1$, with $\max_{0\le t\le 1}\kappa(g_t)<8\pi+\epsilon$ for any $\epsilon>0$.
\begin{figure}[htbp]
\begin{center}
\includegraphics[angle=0, width=10cm]{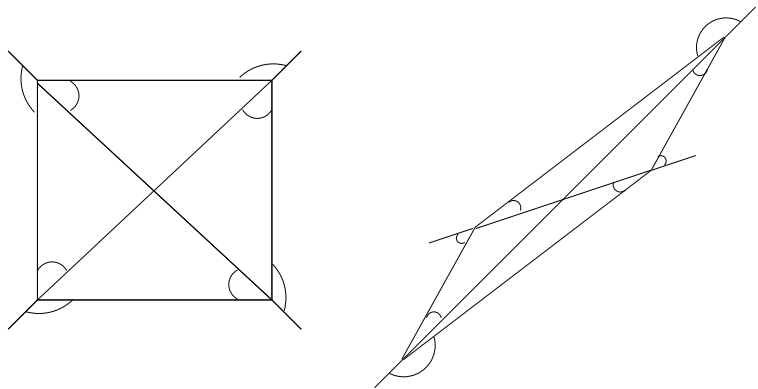}
\end{center}
\caption{Deformation of the degenerate middle stage}
\label{fig:inout3}
\end{figure}
To finish the proof we note that \fullref{lmacountmax} implies that $g_0$ and $g_1$ can be taken to have total curvature arbitrarily close to $4\pi$. It will thus be sufficient to connect these immersions with standard spheres keeping $\kappa$ smaller than $8\pi$. Regarding the elongated immersions with two extrema in the long direction as a family of immersed plane curves as in the proof of \fullref{lmatwoextr} it is not hard to show that it is possible to connect these to standard sphere keeping $\kappa$ close to $4\pi$ by elongating the entire deformation in the direction where the spheres have only two extrema until we arrive at a sphere immersion which is the boundary of a convex body. The latter can then be connected to a standard sphere through boundaries of convex bodies.
\end{proof}

\section{General properties of local minima}\label{sec:7}
In this section we discuss local minima of the total curvature functional in general dimensions using some tools from contact geometry.

\subsection{Basic contact notions}\label{ssec:7A}
A {\em contact structure} on an orientable $(2n+1)$--manifold $N$ is a completely non-integrable field of tangent hyperplanes $\xi\subset TN$. That is, a field of hyperplanes given as $\xi=\krn{\alpha}$, where the non-vanishing $1$--form $\alpha$ (the {\em contact form}) is such that the $(2n+1)$--form $\alpha\wedge(d\alpha)^n$ is a volume form on $N$. Note that if $\alpha$ is a contact form then $d\alpha|\xi$ is a symplectic form. A diffeomorphism of contact manifolds $(N,\xi)\to (N',\xi')$ is called a {\em contactomorphism} if it maps $\xi$ to $\xi'$. An immersion of an $n$--manifold $f\co L\to N$ is called {\em Legendrian} if $df_p(T_p M)\subset \xi_{f(p)}$ for all $p\in L$.

\subsection{The Legendrian lift of an immersion}\label{ssec:7B}
Let $M$ be an $n$--manifold, let $T^\ast M$ be its cotangent bundle, and let $\pi\co T^\ast M\to M$ be the projection. Define the $1$--form $\theta_M$ on $T^\ast M$
$$
\theta_M(p)V=p(d\pi V),\quad p\in T^\ast M,\,\, V\in T_p(T^\ast M).
$$
The unit cotangent bundle $UT^\ast\R^{n+1}$ of $\R^{n+1}$ carries a natural contact structure. Consider $UT^\ast\R^{n+1}$ as a subset of $T^\ast\R^{n+1}$, then the restriction $\alpha$ of the $1$--form $\theta_{\R^{n+1}}$ to $U T^\ast\R^{n+1}$ is a contact form. Also the $1$--jet space $J^1(\S^n,\R)=T^\ast\S^n\times\R$ of $\S^n$ carries a natural contact structure with contact $1$--form
$$
\beta = dz-\theta_{\S^n},
$$
where $z$ is a linear coordinate on the $\R$--factor.

If $v\in\R^{n+1}$ then let $v^\ast$ be the linear form on $\R^{n+1}$ given by $\la v,\cdot\ra$, and if $p$ is a linear form let $p^\ast$ be the vector such that $p=\la p^\ast,\cdot\ra$. We write elements in $T^\ast\R^{n+1}$ as $(x,p)$, where $x\in\R^{n+1}$ and where $p$ is a linear functional on $\R^{n+1}$. Likewise, elements in $T^\ast_v\S^n$ will be written as $(v,p)$ where $v\in\S^{n}\subset\R^{n+1}$ and where $p$ is a linear functional on $\R^{n+1}$ such that $p(v)=0$. Note that the diffeomorphism $\Phi\co U T^\ast\R^{n+1}\to J^1(\S^n)$
$$
\Phi(x, p)=\Bigl(p^\ast,\,x^\ast - p(x)p,\, p(x)\Bigr)
$$
is a contactomorphism.

If $f\co M\to\R^{n+1}$ is an immersion of an oriented $n$--manifold we define its {\em Legendrian lift} $Lf\co M\to UT^\ast\R^{n+1}$ by
\begin{equation}\label{eqn:leglift}
Lf(p)=\Bigl(f(p),\nu(p)^\ast\Bigr),
\end{equation}
where $\nu(p)$ is the positive normal of $f$. Using $\Phi$ we may regard $Lf$ as a map into $J^1(\S^n)$ as well. Note that $Lf$ is a Legendrian immersion.

Let $\pi_1\co UT^\ast\R^{n+1}\to\R^{n+1}$ and $\pi_2\co UT^\ast\R^{n+1}\to \S^{n}$ be the natural projections. If $\Gamma\co M\to UT^\ast\R^{n+1}$ is a Legendrian immersion such that $\pi_1\circ\Gamma$ is an immersion then there exists an orientation on $M$ such that $\pi_2\circ\Gamma^\ast$ is the Gauss map of $\pi_1\circ\Gamma$.

Finally, let $\phi\co \S^{n}\to \S^{n}$ be a diffeomorphism then $\phi$ has a natural lift to a contactomorphism $\Phi\co J^1(\S^n)\to J^1(\S^n)$,
\begin{equation}\label{eqn:contdiff}
\Phi(x,p,z)=(\phi(x),p\circ [d\phi^{-1}],z).
\end{equation}
Note that $\Phi$ depends continuously on $\phi$.

\subsection{Shrinking volumes of maps and of Gauss maps}\label{ssec:7C}
Let $M$ be a closed $n$--manifold and let $f\co M\to\S^n$ be any smooth map. By Sard's theorem the critical values of $f$ form a subset of $\S^n$ of measure zero. Assume that there exists a point $p\in M$ such that $\rank(df_p)=n$. Then the volume $\vol(f)$ of the map $f$ satisfies $\vol(f)>0$. If $p\in\S^n$ is a regular value we say that the {\em absolute multiplicity} of $f$ at $p$ is the (finite by compactness) number of points in $f^{-1}(p)$. We prove two lemmas about decreasing volumes of maps.
\begin{iilma}\label{lmamult}
Assume that $f\co M\to\S^n$ has two regular values (a non-value is a regular value of absolute multiplicity $0$) of different absolute multiplicities. Then there exists a $1$--parameter family of diffeomorphisms $\phi_t\co \S^n\to \S^n$, $0\le t\le 1$, with $\phi_0=\id$ and such that
$$
\vol(\phi_t\circ f)<\vol(f), \text{ for }0<t\le 1.
$$
\end{iilma}

\begin{proof}
Consider first the $1$--dimensional case. Let $p$ and $q$ be points in $\S^1$ of absolute multiplicity $m$ and $n$ respectively with $m<n$. Let $A\subset\S^1$ be the positively oriented arc connecting $p$ to $q$. Since the subspace of regular values in $\S^1$ is open and since the absolute multiplicity is locally constant on this subspace, there exists $\delta>0$ such that in a $\delta$--neighborhood of $\pa A$ in $A$ the multiplicity is constant. Let $B(x,r)$ denote the $r$--ball around $x$ and choose a smooth function $\phi\co\S^1\to [0,1]$ with the following properties
\begin{itemize}
\item $\phi=0$ on a $\frac18\delta$--neighborhood of $\S^1-A$,
\item $\phi=1$ on $A-\left( B(p,\frac12\delta)\cup B(q,\frac12\delta)\right)$,
\item $\frac{d\phi}{ds}\ge 0$ on $B(p,\delta)$ and
$\frac{d\phi}{ds}\le 0$ on $B(q,\delta)$.
\end{itemize}
For small $\theta>0$ consider the $1$--parameter family of diffeomorphisms $\phi_t\co\S^1\to\S^1$
$$
\phi_t(x)=e^{i\theta\phi(x)t}x.
$$
If $\theta>0$ is sufficiently small then, for $t>0$, $\phi_t$ strictly decreases the length of the region where the multiplicity equals $n$ and strictly increases the length of the region where the multiplicity equals $m$. The lengths of all regions of other multiplicities are left unchanged. The lemma follows in the $1$--dimensional case.

In the higher dimensional case. Let $p, q\in\S^n$ be regular values of multiplicities $m$ and $n$, respectively with $m<n$. Connect $p$ to $q$ by an oriented  great circle arc $A\subset\S^1\subset\S^n$.  Again, for sufficiently small $\delta>0$, the absolute multiplicity is constant on a $\delta$--neighborhood of $\pa A$. Consider a tubular neighborhood $T\approx \S^1\times D(\eta)$, with fiber disks $D(\eta)$ of radii $0<\eta<\delta$, of the great circle $\S^1$. Let $s$ be a coordinate on $\S^1$ and let $\phi(s)$ be a function as above. Let $g\co [0,\eta]\to [0,1]$ be a non-increasing smooth function with $g(0)=1$ and $g=0$ on $[\frac12\eta,\eta]$. Let $\psi\co D(\eta)\to[0,1]$ be the function $\psi(\xi)=g(|\xi|)$. Let $(s,\xi)\in\S^1\times D(\eta)$ be coordinates on $T$. For small $\theta>0$ define the $1$--parameter family of diffeomorphism $\phi_t\co\S^n\to\S^n$ as $\phi_t=\id$ on $\S^n-T$ and
$$
\phi_t(s,\xi)=\Bigl(e^{i\theta\phi(s)\psi(\xi)t}s,\xi\Bigr),
$$
on $T$. Since the metric of $\S^n$ has the form
$$
ds^2+f(|\xi|)d\xi^2
$$
in $T$, it follows that, for $t>0$, $\phi_t$ strictly decreases the volume of the set of regular values of absolute multiplicity $n$, strictly increases the volume of the set of regular values of absolute multiplicity $m$, and leaves the volumes of all regions of other multiplicities invariant. The lemma follows.
\end{proof}

\begin{iilma}\label{lmafold}
Assume that $f\co M\to\S^n$ is a smooth map and that there are local coordinates $x=(x',x'')\in\R\times\R^{n-1}$ around $p\in M$ and coordinates $y=(y',y'')\in\R\times\R^{n-1}$ around $f(p)\in\S^n$ such that
$$
f(x',x'')=\left((x')^{2k},x''\right), \text{ for some }k>0.
$$
Then there is a $1$--parameter family of diffeomorphisms $\phi_t\co\S^n\to\S^n$, $0\le t\le 1$, supported in a small neighborhood of $f(p)$, and a cut-off function $\psi(x)$ supported in a neighborhood of $x$ such that
$$
\vol(\phi_{t\psi}\circ f)<\vol(f),\quad \text{for } 0<t\le 1.
$$
\end{iilma}

\begin{proof}
The lemma can be proved by an obvious modification of the proof of \fullref{lmamult}.
\end{proof}

\begin{iiprp}\label{prpGauss1}
Let $f\co M\to\R^{n+1}$ be an immersion which is a local minimum of $\kappa$ of normal degree $d$. Then $\kappa(f)=\vol(\S^n)(|d|+2k)$, for some integer $k\ge 0$.
\end{iiprp}

\begin{proof}
Note that the Gauss map $\nu$ of $f$ is smooth and that any regular value has algebraic multiplicity $d$, and therefore absolute multiplicity $|d|+2k$ for some $k\ge 0$. If the volume of the Gauss map is different from $\vol(\S^n)(|d|+2m)$ for all integers $m$ then it follows that there are regular values of $\nu$ of different multiplicities. Let $\phi_t\co\S^n\to\S^n$, $0\le t\le 1$, be a family of diffeomorphisms as in \fullref{lmamult} which shrinks $\vol(\nu)$. This family induces a family of contactomorphisms $\Phi_t\co J^1(\S^n)\to J^1(\S^n)$, see \eqref{eqn:contdiff}, and hence of $UT^\ast\R^{n+1}$. Since the set of Legendrian immersions $\Gamma\co M\to UT^\ast\R^{n+1}$ such that $\pi_1\circ\Gamma$ is an immersion is open, it follows that for $\epsilon>0$, sufficiently small, $\pi_1\circ\Phi_t\circ f$, $0\le t\le \epsilon$, is a regular homotopy of $f$ shrinking the total  curvature. The proposition follows.
\end{proof}

\begin{iiprp}\label{prpGauss2}
If the Gauss map of $f\co M\to\R^{n+1}$ has the form of the map in \fullref{lmafold} then there exists a regular homotopy $f_t$, $0\le t\le 1$, with $f_0=f$ and such that
$$
\frac{d}{dt}\kappa(f_t)|_{t=0}<0.
$$
\end{iiprp}

\begin{proof}
Similar to the proof of \fullref{prpGauss1}.
\end{proof}

\section{RIT immersions and curvature generic local minima}\label{sec:8}
In this section we discuss relatively isotopy tight (RIT) immersions. These are closely related to local minima of the total curvature functional. We prove a non-injectivity result for RIT immersions and apply it to demonstrate \fullref{thmgenlocmin}.

\subsection{Relatively isotopy tight immersions}\label{ssec:7D}
Let $M$ be a closed surface and recall that if $f\co M\to\R^3$ and $v\in\S^2$ then $f_v$ denotes the height function in direction $v$ composed with $f$.
\begin{iidfn}\label{dfnRIT}
An immersion $f\co M\to\R^3$ is {\em relatively isotopy tight (RIT)} if
\begin{equation*}
\kappa(f)=2\pi\max_{v\in\S^2}\mu(f_v)=
2\pi\min_{v\in\S^2}\mu(f_v).
\end{equation*}
and
\begin{equation*}
\min_{v\in\S^2}\mu(g_v)\ge\frac{\kappa(f)}{2\pi},
\end{equation*}
for every immersion $g$ in some neighborhood of $f$.
\end{iidfn}

This is Definition A, in Kuiper and Meeks \cite{iiKuMe}, where the following structure theorem appears as Fundamental Lemma 2. (For simplicity of formulation we state it for embeddings, the obvious analog for immersions holds as well.)

\begin{iithmKM}
If $f\co M\to\R^3$ is an RIT $C^2$--embedding (immersion) then there exists an integer $k$ and convex surfaces $\pa B_0,\dots,\pa B_k$ in $\R^3$ such that the set $M_{K>0}=K^{-1}(0,\infty)$ satisfies
$$
M_{K>0}=\bigcup_{j=0}^k (\pa B_j)_{K>0}.
$$
The unique principal component $M_j^+$ of $\pa B_j\cap M$ contains $(\pa B_j)_{K>0}$ and is obtained from $\pa B_j$ by deleting disjoint plane convex disks in $\pa B_j$. The boundary of such a disk is called a {\em top-circle}. The plane $\Pi(\gamma)$ of a top-circle $\gamma$ supports some neighborhood $U\subset M$ of $\Pi(\gamma)\cap M_j^+$. Moreover,
$$
\kappa(f)=4\pi(1+g(M)+2k).
$$
All these properties together are sufficient for an embedding (immersion) to be RIT.
\end{iithmKM}

\subsection{Non-injectivity}\label{ssec:7E}
In \cite{iiKuMe}, the following question concerning RIT immersions of $2$--spheres is posed as Problem 5: {\em Is there an RIT embedding of the $2$--sphere with $\kappa=4\pi+8\pi k$ for some $k>0$?}, and it is shown that there are no such embeddings for $k=1$. In this subsection we prove the following result which give some partial information about this question (and which, as we shall see in \fullref{ssec:7F}, also leads to a proof of \fullref{thmgenlocmin}).

\begin{iithm}\label{thmRIT}
Let $f\co S^2\to\R^3$ be an RIT immersion with the following properties.
\begin{itemize}
\item[{\rm (a)}] At least one component of $K^{-1}((-\infty,0])$ which intersects $K^{-1}((-\infty,0))$ is an annulus with one boundary component bounding a disk $D\subset S^2$ such that $K\ge 0$ on $D$.
\item[{\rm (b)}]
If $\nu$ is a non-vanishing normal vector field of $f$ and if $p,q\in S^2$ are any two points with $K(p)>0$ and $K(q)>0$ then $\la \nu(p),H(p)\ra\cdot\la \nu(q),H(q)\ra>0$, where $H(p)$ is the mean curvature vector of $f$ at $p$. (In other words the signs of $\la \nu(p),H(p)\ra$ agree for all $p\in S^2$ with $K(p)>0$.)
\end{itemize}
Then $f$ is non-injective. It follows in particular that every RIT immersion such that all components of $K^{-1}((-\infty,0])$ are annuli is non-injective.
\end{iithm}

\fullref{thmRIT} will be proved using two lemmas which we present next. It should be viewed in the light of the structure theorem for RIT immersions which has the following consequences. Let $f\co S^2\to\R^3$ be an RIT immersion. Then the structure theorem in \fullref{ssec:7D} gives a subdivision of $S^2$ into planar surfaces of two kinds. The first kind is the planar surfaces $\Omega_j^+$ which map in a one to one fashion onto the convex bodies $\pa B_j$ with planar convex disks removed. We write $\Omega^+$ for the union $\cup_j \Omega_j^+$ over all components. The second kind of planar surface is a component of $S^2-\Omega^+$. We denote such components $\Omega_k^-$ and their union $\Omega^-=\cup_k\Omega_k^-$. Then each $\Omega_k^-$ is a sphere with $r\ge 2$ disks removed. Our first lemma gives a non-embeddedness condition for such components which are spheres with two disks removed.

Let $f\co S^1\times[0,1] \to\R^3$ be an immersion with the following properties.
\begin{itemize}
\item[(i)] The map $f|S^1\times\{a\}$ is a planar curve bounding a convex region $D_a$ in a plane $\Pi_a$, $a=0,1$,
\item[(ii)] for some neighborhood $N_a$ of $f^{-1}(\Pi_a)$, $f(N_a)$ is supported by $\Pi_a$, and adjoining to $f$ a small collar in the plane $\Pi_a$, outside $D_a$, we get a $C^1$--immersion, $a=0,1$, and
\item[(iii)] $\kappa(f)=4\pi$.
\end{itemize}
The Gauss--Bonnet theorem then implies that $K\le 0$ everywhere in $S^1\times[0,1]$. It follows from this that $f(S^1\times[0,1])$ is contained in the convex hull of $\pa D_0$ and $\pa D_1$ (if not, it is straightforward to find a point $p\in  f(S^1\times[0,1])$ where $K>0$).

\begin{iilma}\label{lmacylinj}
If $D_0\cap D_1\ne \emptyset$, and if $\pa D_0$ and $\pa D_1$ are unlinked then $f$ is non-injective.
\end{iilma}

\begin{proof}
We have $D_0\cap D_1\ne\emptyset$ and $K\le 0$ in $S^1\times[0,1]$. If $\pa D_0\cap\pa D_1\ne\emptyset$ then $f$ is not injective. Thus assume that $\pa D_0\cap \pa D_1=\emptyset$. The linking condition then implies that the boundary of one of the disks does not intersect the other. Assume for definiteness that $D_0\cap \pa D_1=\emptyset$.

We claim that $f(S^1\times[0,1])\cap\inr(D_0)\ne\emptyset$, where $\inr(X)$ denotes the interior of $X$. Suppose that this is not the case. Then we complete $f(S^1\times[0,1])$ to an embedding $g\co D^2\to\R^3$ of the $2$--disk $D^2$ as follows. Add a small collar $C_0(\eta)=\pa D_0\times[0,\eta]$ along $f(S^1\times\{0\})$ in the direction of the normal $\nu_0$ of the plane $\Pi_0$ which points into the half space of $\R^3$ not containing $f(N_0)$. Also add a copy of $D_0$, $D_0(\eta)$ shifted $\eta$ units along this normal, and finally add a large annular region $A_1$ in $\Pi_1$ bounded by $\pa D_1$ and the boundary of a large disk containing $D_1$. This is the image of the piecewise smooth embedding $g\co D^2\to\R^3$. (This is an embedding since $f(S^1\times[0,1])$ lies inside the convex hull of its boundary.)

Let $h_1$ denote the height function in the direction perpendicular to $\Pi_1$ normalized so that $h_1=0$ on $\Pi_1$. We can further complete the embedding $g$ of the disk to embeddings of a sphere in two different ways: add the lower hemisphere ($h_1\le 0$) of the sphere containing $\pa A_1-\pa D_1$ as a great circle, or add the upper hemisphere ($h_1\ge 0$) of this sphere. We call the former embedded sphere $G_l$ and the latter $G_u$.

Then $G_u$ and $G_l$ bound balls in $\R^3$. Noting that the global maximum on $G_l$ of $h_1$ and the global minimum on $G_u$ of $h_1$ lies in $D_0(\eta)\cup C_0(\eta)$ we find that the outward coorientations along $D_0(\eta)$ of both $G_u$ and of $G_l$ must point along $\nu_0$. This local coorientation however determines the outward coorientations of $G_u$ and $G_l$ along $A_1$. In particular, the outward coorientations of $G_u$ and $G_l$ along $A_1$ must agree. This, however, contradicts the fact that the global minimum of $h_1$ on $G_l$ and the global maximum of $h_1$ on $G_u$ are attained on the added half-spheres since their outward coorientations induce opposite coorientations on $A_1$. We conclude that $f(S^1\times I)\cap\inr(D_0)\ne \emptyset$.

After an arbitrarily small translation of $D_0$ along the normal line to $D_0$ we may assume that $f$ intersects $\inr(D_0)$ transversely. If some component of this intersection bounds a disk in $S^1\times[0,1]$, then this disk has an extremum in a direction perpendicular to $D_0$ and thus has a non-degenerate extremum in some direction arbitrarily close to this one. Such an extremum contradicts $K\le 0$. On the other hand, if no component of the intersection bounds a disk in $S^1\times[0,1]$ then there exists a component which together with $\pa D_0$ bounds a cylinder which must have at least one extremum in the direction perpendicular to $D_0$, which again contradicts $K\le 0$. We conclude that $f$ cannot be an embedding.
\end{proof}

Let $f\co S^1\times[0,1] \to\R^3$ be an immersion with the properties (i)--(iii) (stated above \fullref{lmacylinj}) which is also an embedding. Let $\nu_0$ and $\nu_1$ be the unit normal vectors of $\Pi_0$ and $\Pi_1$ which point into the half spaces which do not contain $f(N_0)$ and $f(N_1)$, respectively. Moreover, assume that $\pa D_0$ and $\pa D_1$ are unlinked.

\begin{iilma}\label{lmacylcoor}
If $\nu$ is a normal vector field of $f$ then for any $p\in S^1\times\{0\}$ and any $q\in S^1\times\{1\}$ $\la \nu_0,\nu(p)\ra>0$ if and only if $\la \nu_1,\nu(q)\ra>0$.
\end{iilma}

\begin{proof}
For small $\eta>0$, let $C_j(\eta)=\pa D_j\times[0,\eta]$ be the collar on $\pa D_j$ in direction $\nu_j$ and let $D_j(\eta)$ be $D_j$ shifted $\eta$ units along $\nu_j$. It follows from \fullref{lmacylinj} that for $\eta>0$ small enough
$$
D_0(\eta)\cup C_0(\eta)\cup f(S^1\times I)\cup C_1(\eta)\cup D_1(\eta)
$$
is an embedded sphere. Considering the height functions in directions perpendicular to $\Pi_0$ and $\Pi_1$ we find that the outward coorientation of this sphere is $\nu_0$ along $D_0(\eta)$ and $\nu_1$ along $D_1(\eta)$. The lemma follows from the continuity of the coorientation.
\end{proof}

\begin{proof}[Proof of \fullref{thmRIT}]
We argue by contradiction: assume we can find an RIT embedding $f$ with properties as described. Note first that if any two top-circles in the decomposition of $f$ are linked then $f$ cannot be an embedding. Assume thus these circles are pairwise unlinked.

Consider the image $f(D)$ of the disk $D$ where the curvature is non-negative, $D\subset\Omega_+$. Let $B$ denote the convex hull of $f(D)$. Let $A\subset \Omega_-$ denote the non-positively curved cylinder in the decomposition of $f$ such that $\pa B\cap \pa A\ne \emptyset$. Write $\pa A= \Sigma_1\sqcup \Sigma_2$ and let $\Sigma_1=\pa B\cap \pa A$. \fullref{lmacylinj} implies that the convex planar disk $\Delta_2$ bounded by $\Sigma_2$ must either lie entirely outside $B$ or entirely inside $B$.

If $\Delta_2$ lies outside $B$ then $A\cap \inr(\Delta_1)\ne \emptyset$, where $\Delta_1$ is the convex planar disk bounded by $\Sigma_1$. (To see this, note that points in $A$ near $\Sigma_1$ lies inside $B$.) The argument used in the proof of \fullref{lmacylinj} then shows there exists a point where $K>0$ in $A$. But $K\le 0$ in $A$, hence $\Delta_2$ lies inside $B$.

As in the proof of \fullref{lmacylcoor} let, for small $\eta>0$, $C_j(\eta)=\pa\Sigma_j\times[0,\eta]$, $j=1,2$, be a collar in the direction of the normal of the plane of $\Delta_j$ which points away from $A$ and let $\Delta_j(\eta)$ be $\Delta_j$ shifted $\eta$ units along this normal. Note that $f(D)\cup A\cup C_2(\eta)\cup \Delta_2(\eta)$ is an embedded $2$--sphere which subdivides $\R^3$ into two connected components. One of these components is a subset of $\inr(B)$, we call that component $X(\eta)$.

Let $\Pi_1$ and $\Pi_2$ be the planes which contain $\Delta_1$ and $\Delta_2$, respectively.
Consider an embedding $k\co D^2\to \R^3$  with the following properties:
\begin{itemize}
\item[(i)] $k(\pa D^2)=\Sigma_2$,
\item[(ii)] there exists some neighborhood $U$ of $k^{-1}(\Pi_2)$ and some
neighborhood $V$ of $\Pi_1\cap A$ in $A$ such that
$k(U)\cup V$ is a smooth surface supported by $\Pi_2$, and
\item[(iii)] $k(\inr (D^2))\cap (f(D)\cup A\cup C_2(\eta)\cup\Delta_2(\eta))=\emptyset$.
\end{itemize}
Clearly, such embeddings $k$ satisfies $k(D^2)\subset \overline{X(\eta)}$.

Let $E=S^2-(D\cup f^{-1}(A))$. Then $f\co E\to\R^3$ is a disk embedding which satisfies (i) and (ii). Moreover, for almost all $\eta>0$, $f(E)$ intersects $\Delta_2(\eta)$ transversely. We construct an embedding $g\co D^2\to\R^3$ fulfilling also (iii) as follows.

By transversality $f(E)\cap \inr(\Delta_2(\eta))$ is a finite collection of circles. Consider an innermost circle $\gamma$ in this intersection and its preimage $\Gamma$ in $E$. Note that $\Gamma$ divides $E$ into a disk $E_d$ and an annulus $E_a$. Define the map $g'\co D^2\to\R^3$ by cutting $f(E)$ along $\gamma$ and replacing $f(E_d)$ with the disk bounded by $\gamma$ in $\Delta_2(\eta)$ shifted slightly in the direction of $f(E_a)$. By induction we remove all intersection circles. Let $g\co D^2\to \overline{X(\eta)}$ denote the embedding constructed in this way.

Let $h_2$ be the height function in the direction perpendicular to $\Pi_2$ such that $h_2(\Pi_2)=0$, and such that $h_2\circ g$ is positive near the boundary of $D^2$. Then $h_2\circ g$ has a global minimum on $D^2$. Assume first that the value of $h_2\circ g$ at its global minimum equals $0$.

Note that $\Delta_1(\eta)\cup C_1(\eta)\cup A\cup C_2(\eta)\cup\Delta_2(\eta)$ is an embedded sphere which bounds a ball. Let $\alpha(\eta)$ be an arc inside this ball connecting $\Delta_1(\eta)$ to $\Delta_2(\eta)$ and meeting the boundary transversely. Add to $\alpha(\eta)$ two half-rays perpendicular to $\Delta_1(\eta)$ and $\Delta_2(\eta)$, respectively, to get an embedded curve $\beta(\eta)$. Note that $g(\pa D^2)$ and $\beta(\eta)$ have linking number one. Therefore $g(D^2)$ must intersect $\beta(\eta)$. In fact it must intersect $\alpha(\eta)$ since the global minimum of $h_2\circ g$ equals $0$. This however contradicts $k(D^2)\subset \overline{X(\eta)}$ for all $\eta>0$ since $\alpha(\eta)\in\R^3-\inr(X(\eta))$.

It follows that the value at the global minimum of $h_2\circ g$ is smaller than $0$. This implies in particular that the global minimum on $g(D^2)$ of a height function $h$ very close to $h_2$ is attained at some point of $q$ of $g(D)$ with $K(q)>0$. Since the part of the image of $g$ which is not in the image of $f$ can be taken arbitrarily close to $\{h_2=0\}$ this non-degenerate global minimum $q$ must be a point in the image of $f$. Consider the embedded sphere $Y=f(D)\cup A\cup g(D^2)$. Since $q$ is a global minimum in $g(D^2)$ of a height function $h$ arbitrarily close to $h_2$ there exists a path in the closure of the bounded component of $\R^3-Y$ which connects $q$ to a point $q'\in f(D)$ such that $h(q')<h(q)$ and such that $K(q')>0$. It follows that the outward coorientation $\nu$ of $Y$ satisfies
\begin{equation}\label{eqn:wrong}
\la\nu(q),H(q)\ra=-\la\nu(q'),H(q')\ra.
\end{equation}
However, it is clear from the above construction that the coorientation $\nu$ of $Y$ agrees (up to an {\em over all} sign) with the coorientation of $f$ for points in $Y$ which are in the image of $f$. Thus, \eqref{eqn:wrong} contradicts our coorientation assumption on $f$ and we conclude that $f$ cannot be an embedding.

The last statement of the theorem is a direct consequence of the above and \fullref{lmacylcoor}
\end{proof}

\subsection{Curvature generic local minima}\label{ssec:7F}
Let $f\co S^2\to\R^3$ be an immersion. Let $K$ denote the curvature function of the metric induced by $f$. We say that $f$ is {\em curvature generic} if $0$ is a regular value of $f$.

\begin{iilma}\label{lmagenlocmin}
Let $f\co S^2\to\R^3$ be a curvature generic local minimum of $\kappa$. Then $f$ is an RIT immersion such that any component $\Omega_j^-$ of $K^{-1}(-\infty,0]$, see \fullref{ssec:7E} for notation, is a two-holed sphere.
\end{iilma}

\begin{proof}
Let $\gamma\subset S^2$ be a component of $K^{-1}(0)$. Then $\gamma$ is an embedded circle. Let $s$ be a parameter along $\gamma$. We show that
$$
\frac{d}{ds}\nu(s)=0,
$$
where $\nu$ is the Gauss map of $f$. Assume not, then there exists an open arc $\alpha$ in $\gamma$ where $\frac{d}{ds}\nu\ne 0$. Since $f$ is curvature generic there exists coordinates $(u,v)\in(-\delta,\delta)^2=W$ around some point $p\in \alpha$ with $K^{-1}(0)=\{v=0\}$.  Since $\frac{d}{ds}\nu(p)\ne 0$, the second fundamental form of $f$ is non-trivial at $p$ and hence has a non-zero eigenvalue. But since $K(p)=0$ at least one of the eigenvalues of the second fundamental equal zero. It follows that if $W$ is chosen small enough then the eigenvalues of the second fundamental form of $f$ are distinct throughout $W$. Let $e_1,e_2$ be the eigenvector fields of the second fundamental form along $W$, with corresponding eigenvalues $\lambda_1,\lambda_2$. Then $K=\lambda_1\lambda_2=0$ along $\{v=0\}$. Since $0$ is a regular value of $K$ we find
$$
\frac{\pa\lambda_1\lambda_2}{\pa v}|_{(u,0)}=
\frac{\pa\lambda_1}{\pa v}|_{(u,0)}\lambda_2(u,0)
+\lambda_1(u,0)\frac{\pa\lambda_2}{\pa v}|_{(u,0)}\ne 0.
$$
This implies that at most one of $\lambda_1,\lambda_2$ vanishes along $\{v=0\}$. Assume it is $\lambda_1$. Since $\frac{\pa}{\pa u}\nu(u,0)\ne 0$ we see that $e_1$ and $\pa_u$ are linearly independent. Moreover $\frac{\pa\lambda_1}{\pa v}(u,0)\ne 0$. Using the flow of $e_1$ we construct a coordinate system $(u,x)$ in a neighborhood of $\{v=0\}$. We have
$$
\nu(u,x)=\nu(u,0)+\frac{\pa}{\pa x}\nu(u,0)x
+\frac12\frac{\pa^2}{\pa x^2}\nu(u,0)x^2+\Ordo(x^3).
$$
Since $\frac{\pa}{\pa x}\nu(u,x)=\lambda_1e_1$ and $\frac{\pa^2}{\pa x^2}\nu(u,0)=\frac{\pa\lambda_1}{\pa x}e_1$,
$$
\nu(u,x)=\nu(u,0)+
+\frac12\frac{\pa\lambda}{\pa x}e_1(u,0) x^2+\Ordo(x^3).
$$
Thus $\nu$ has a fold singularity along $\{v=0\}$. By \fullref{prpGauss2}, this contradicts $f$ being a local minimum and we conclude that $\nu(s)=\nu_\gamma$ is constant along $\gamma$. Since $\la \dot\gamma,\nu_\gamma\ra=0$ we find that $\gamma$ lies in a plane orthogonal to $\nu_\gamma$.

We subdivide the components of $K^{-1}(0)$ into levels as follows. Any component $\gamma$ subdivides the sphere into two disks. If one of these disks do not contain other components of $\Omega_0=K^{-1}(0)$ we say that $\gamma$ is a curve of level $0$. In general we make the following inductive definition. Let $\Omega_j$ be the subset of $K^{-1}(0)$ obtained by removing from it all components of level smaller than $j$. Define a component of $K^{-1}(0)$ to have level $j$ if one of the two disks into which it divides $S^2$ does not contain any component in $\Omega_j$.

Let $\gamma$ be a $0$--level component of $K^{-1}(0)$. By the above
$f(\gamma)$ is planar curve in a plane with normal $\nu_\gamma$. We
claim that $f(\gamma)$ must be convex. Consider an arc $A$ of
$\gamma$. The union of one piece of the surface near $A$ (the
positively curved one) and a planar region with boundary $A$ must be
convex. It is clear that this piece lies in the direction of the
curvature vector of $A$. Hence, since $\gamma$ separates negative
curvature from positive (because $0$ is a regular value of $K$), it
must be locally convex and the plane in which it lies is a local
support plane of $f$. Moreover, restricted to any curve sufficiently
near $\gamma$ the image of the Gauss map is a curve wrapping around
$\nu_\gamma\in \S^2$ $n$ times. Filling this curve with a disk with
singularity in the middle we get a branched cover of the sphere by a
sphere with one singular point of multiplicity $n$. The
Riemann--Hurwitz formula gives
$$
\chi(S^2)=n\chi(S^2)-(n-1)=n+1,
$$
where $\chi$ denotes the Euler characteristic. Hence $n=1$ and $f(\gamma)$ is convex.

Assume inductively that the images of all curves of level $j-1$ are planar and convex. Consider a level $j$ curve $\pa R_+$ which bounds a region $R\subset S^2$ and such that all other boundary components $\pa R_-$ of $R$ are curves of level $j-1$. Consider the Gauss map restricted to $R$. Our inductive assumption shows that the map $\nu|(R-\pa R_+)$ gives an immersion from the disk $\bar R$ obtained from $R$ by filling each component of $\pa R_-$ with a disk. Arguing as above we find that $f(\pa R_+)$ is convex as well and that $\kappa(R)=\area(\nu(\bar R))=4\pi$. Thus by the Gauss--Bonnet theorem, if $h$ denotes the number of holes in the negatively curved sphere with holes
$$
\int_R K\,dA=-4\pi=2\pi(2-h) - 2\pi h=4\pi-4\pi h.
$$
Hence $h=2$ and it follows that each negatively curved component is a cylinder. (One may see that each component of $K^{-1}((-\infty,0))$ is a cylinder also by using $K<0$ and the line fields arising as asymptotic directions.)
\end{proof}

\begin{proof}[Proof of \fullref{thmgenlocmin}]
It follows from \fullref{lmagenlocmin}, that any curvature generic local minimum of $\kappa$ is an RIT immersion which satisfies the conditions of \fullref{thmRIT}. The theorem then follows from \fullref{thmRIT}.
\end{proof}

\end{document}